\documentclass[francais]{cedram-aif}

\usepackage[english,frenchb]{babel}
\usepackage[latin1]{inputenc}
 \usepackage[T1]{fontenc}
 \usepackage{lmodern}

% ...

% ...

\newcommand{\dem}{\textbf{Proof: }}
\newcommand{\rem}{\textbf{Remark: }}

\newtheorem{ass}{Assumption}

\equalenv{cor}{coro}
\equalenv{lem}{lemm}
\equalenv{definition}{defi}

\title
[Normal form of vector fields with an invariant torus]
{Normal form of holomorphic vector fields with an invariant torus under Brjuno's A condition}

\author{\firstname{Claire} \lastname{Chavaudret}}

\address{Universit\' e de Nice Sophia-Antipolis\\ 
Laboratoire J.A.Dieudonn\' e\\
Parc Valrose\\
06108 Nice Cedex 02 (France)}

\email{claire.chavaudret@unice.fr}

  
\altkeywords{normal forms, invariant torus, Brjuno condition, small divisors, KAM, resonances}


\begin{document}
\begin{abstract}On prouve l'existence d'une forme normale analytique pour certains champs de vecteurs holomorphes au voisinage d'un point fixe et d'un tore invariant. Apr\`es avoir construit une forme normale formelle, on montre que le champ de vecteurs initial peut \^ etre analytiquement normalis\' e sous deux conditions arithm\' etiques et une condition alg\' ebrique, connues comme les conditions $\gamma,\omega$ et $A$ de Brjuno.
\end{abstract}

\begin{altabstract}This article proves the existence of an analytic normal form for some holomorphic differential systems in the neighborhood of a fixed point and of an invariant torus. Once a formal normal form is constructed, one shows that the initial system with quasilinear part $S$ can be holomorphically conjugated to a normal form, i.e a vector field which commutes with $S$, under two arithmetical conditions known as Brjuno's $\gamma$ and $\omega$ conditions, and an algebraic condition known as Brjuno's $A$-condition, which requires the formal normal form to be proportional to $S$. 
\end{altabstract}

\maketitle

\section{Introduction}

The present article considers the problem of normal forms for autonomous differential systems defined on a product of a torus with a disk, i.e systems of the form

\begin{equation}\label{ptdepart}\dot{X}=f(X,Y); \ \dot{Y}=g(X,Y),\  X\in \mathbb{T}^d, Y\in \mathbb{C}^n
\end{equation}

We also assume that $g(X,0)=0$, $f(X,0)=\omega\in\mathbb{R}^d$ and $\partial_Yg(X,Y)_{\lvert Y=0}=\Lambda$ where $\Lambda$ is a diagonal matrix: thus, the system can be viewed as a perturbation of the integrable system

\begin{equation}\label{unperturbed}\dot{X}=\omega;\ \dot{Y}=\Lambda Y 
\end{equation}

\noindent 
where $X \in \mathbb{T}^d, Y \in \mathbb{C}^n$. 
Such systems can appear when considering the restriction of a larger system on an invariant torus supporting a quasiperiodic motion.
We require the system to be analytic in both variables: thus there exists a complex neighborhood $\mathcal{V}$ of $\mathbb{T}^d$ and a neighborhood $\mathcal{W}$ of the origin in $\mathbb{C}^n$ such that $f$ and $g$ are holomorphic on $\mathcal{V}\times \mathcal{W}$. 

The invariant manifolds of the unperturbed system \eqref{unperturbed} can be easily computed and form a foliation of the phase space. 
What happens after perturbation, however, depends on possible resonances between $\omega$ and the spectrum $\sigma(\Lambda)$ of $\Lambda$. Assuming that $\omega$ and $\sigma(\Lambda)$ are jointly non resonant (in a sense that will be made explicit below), then linearization is possible: as shown in \cite{Au13}, an analytic perturbation of \eqref{unperturbed} is analytically linearizable if $\omega$ and $\sigma(\Lambda)$ are non resonant, and if they satisfy some arithmetical conditions known as Brjuno's $\gamma$ and $\omega$ conditions. 

\bigskip
The problem of holomorphically linearizing a system with non resonant linear part in the vicinity of a fixed point was extensively studied: Siegel (\cite{Siegel}) showed that a diophantine condition on the spectrum of the linear part implies the existence of an analytic linearization. Brjuno (\cite{Bruno}) managed to relax the arithmetical condition on the spectrum. Giorgilli-Marmi (\cite{GM10}) obtained a lower bound for the radius of convergence which is expressed in terms of the Brjuno function of the spectrum.

\bigskip
\noindent
Systems with a fixed point and a resonant linear part have also been studied by Brjuno: 
in \cite{Bruno}, it is proved that for a vector field in $\mathbb{C}^n$ with a fixed point and a resonant linear part, a strong algebraic condition on the formal normal form, known as Brjuno's $A$ condition and the well-known arithmetical "$\omega$ condition" are sufficient in order to have an analytic normalization (i.e an analytic change of variables conjugating the initial vector field to a normal form). Also, the $A$ condition and an arithmetical condition $\bar{\omega}$ which is weaker than $\omega$ are necessary for the analytic normalization.

\noindent In \cite{Stolo00}, other links are given between algebraic conditions on a formal normal form and the holomorphic normalizability of a vector field. 

\noindent A refinement of the normal form theory near a fixed point can be found in Lombardi-Stolovitch \cite{LS10}): if some eigenvalues of $\Lambda$ are zero, then the normal form is defined with respect to a given unperturbed vector field which might be nonlinear. This makes it possible to distinguish the dynamics of distinct vector fields even if they have the same degenerate linear part.


\bigskip
\noindent
In the case under consideration, periodicity with respect to the set of variables $X$ implies a slightly different definition of the normal form: normal forms will be taken in the kernel of $ad_S$, where 

\begin{equation}\label{S}S=\sum_{j=1}^d \omega_j\frac{\partial}{\partial X_j}+\sum_{j'=1}^n \lambda_{j'}Y_{j'}\frac{\partial }{\partial Y_{j'}}
\end{equation}

\noindent i.e $S$ is the unperturbed vector field generating the system \eqref{unperturbed}, where $\lambda_1,\dots, \lambda_n$ are the eigenvalues of the matrix $\Lambda$ and $\omega_j$ are the components of $\omega$. More explicitly, a normal form for \eqref{ptdepart} is a formal vector field $NF(X,Y)$ 
which is a formal Fourier series in $X$ and a formal series in $Y$ and whose Fourier-Taylor development, 

$$NF(X,Y)=\sum_{P\in\mathbb{Z}^d} \sum_{Q\in \mathbb{N}^n} NF_{P,Q} e^{i\langle P,X\rangle } Y^Q$$

\noindent only has non zero coefficients $NF_{P,Q}$ with indices $P,Q$ satisfying 

$$i\sum_{j=1}^d P_j\omega_j +\sum_{j=1}^n Q_j \lambda_j =0$$

\noindent (such couples $(P,Q)$ are the resonances between $\omega$ and $\sigma (\Lambda)$), and 
which is conjugate to 
\eqref{ptdepart} by a formal change of variables (the precise meaning of this formal conjugation will be given in section \ref{Notations}).





\bigskip
The tools used when dealing with vector fields close to an invariant torus are similar to the ones that have been used before, in the study of vector fields with a fixed point. The main difference comes from the presence of a Fourier development which will make the small divisors more complicated to deal with. 
In \cite{Meziani}, Meziani considers the real-valued non resonant case and obtains analytic normalization under a Siegel diophantine condition.


Here, we will combine techniques used by Brjuno, Stolovitch (\cite{Stolo00}) and Aurouet (\cite{Au13}) to prove a conjecture by Brjuno in \cite{Bruno-Book}, mentioned in \cite{Au13}, which is analytic normalization for systems of the form \eqref{ptdepart}, where the linear part might be resonant, under Brjuno's $\gamma,\omega$ and $A$ conditions:

\begin{thm}\label{mainresult}Consider the following autonomous differential system:

\begin{equation}
\label{auton}\dot{X}=\omega+F(X,Y);\ \dot{Y}=\Lambda Y+G(X,Y)
\end{equation}

\noindent where $F$ and $G$ are analytic on a neighborhood $\mathcal{V}$ of $\mathbb{T}^d\times \{0\}$, 
$F(X,0)=0, G(X,0)=\partial_YG(X,0)=0, \Lambda $ is a diagonal $n\times n$ matrix and $\omega\in\mathbb{R}^d$. 
Let $S$ be the quasilinear vector field given by \eqref{S}. Assume 
that 
$S$ satisfies Brjuno's "$\gamma$ and $\omega$ conditions" (see assumption \ref{gamma+omega} below) away from the resonances. If a formal normal form $NF$ of \eqref{auton} is such that there exists a formal series $a(X,Y)$ such that $a(0,0)=1$ and 

\begin{equation}\label{hypsurFN}NF(X,Y)= a(X,Y)S(Y)
\end{equation}

\noindent then there exists $\zeta_0>0$ depending only on $\Lambda,\omega,n,d$ such that if the components of $F$ and $G$ are less than $\zeta_0$ in the analytic norm, then \eqref{auton} is holomorphically normalizable in a complex neighborhood $\mathcal{W}$ of $\mathbb{T}^d\times \{0\}$ and the holomorphic normal form has also the form $b(X,Y)S(Y)$.
\end{thm}

\rem The assumption \eqref{hypsurFN} is usually referred to as Brjuno's $A$ condition. Note that $\mathcal{W}$ is strictly included in $\mathcal{V}$ in general: the loss of analyticity comes from the presence of small divisors.




\bigskip
A first step will be the construction of a normal form by a direct method; although the existence of a formal normal form is not a new result,
one does need to give an explicit construction since the definition of a normal form in the present case might not be completely standard because of Fourier series.
It is then proved that the algebraic assumption \eqref{hypsurFN} on one of the normal forms is preserved by the changes of variables that actually appear in the proof. Then, by a Newton method adapted from Stolovitch and Aurouet, one constructs a converging sequence of holomorphic diffeomorphisms conjugating the system to another one which is normalized up to arbitrary order, which gives the analytic normal form.

\bigskip
\textbf{Acknowledgments:} The author is grateful to Laurent Stolovitch for a useful discussion which helped to improve this article. 

\section{Notations and general definitions}\label{Notations}

\subsection{Topological setting}

Let $r>0$; in the following, $\mathcal{V}_r\subset \mathbb{C}^d$ will denote a complex neighbourhood of the $d$-dimensional torus, of width $r$, and $\mathcal{W}_\delta\subset \mathbb{C}^n$ the ball centered at the origin in $\mathbb{C}^n$, of radius $\delta$.

\bigskip
\noindent
In the following, given a function $f$ on $\mathbb{C}^d/\mathbb{Z}^d\times \mathbb{C}^n$, the notation $f_{P,Q}$ will refer to the complex coefficient appearing next to $e^{i\langle P,X\rangle} Y^Q$ in the Taylor-Fourier development of $f$, that is, the Taylor development near the origin w.r.t. the second variable together with the Fourier development with respect to the first variable (such a development might be purely formal). Also, for all indices $P=(P_1,\dots, P_d)\in\mathbb{Z}^d$ and $ Q=(Q_1,\dots, Q_n)\in\mathbb{N}^n$, one will use the notation $\lvert P\rvert :=\sum_{j=1}^d \lvert P_j\rvert $ (where $\lvert P_j\rvert$ is the absolute value of $P_j$) and $\lvert Q\rvert := \sum_{j'=1}^n Q_{j'}$.
 
 \bigskip
 \noindent 
 The space of scalar-valued functions which are holomorphic on $\mathcal{V}_r\times \mathcal{W}_\delta$ is denoted by $C^\omega_{r,\delta}$. It is provided with the weighted norm
 
 $$\lvert f\rvert_{r,\delta} = \sum_{P\in \mathbb{Z}^d}\sum_{Q\in\mathbb{N}^n} \lvert f_{P,Q}\rvert e^{r\lvert P\rvert} \delta ^{\lvert Q\rvert}$$
 
 \noindent 
 where $f(X,Y)=\sum_{P\in \mathbb{Z}^d}\sum_{Q\in\mathbb{N}^n} f_{P,Q} e^{i\langle X,P\rangle } Y^Q$ is the Taylor-Fourier development of $f$.
 
 \bigskip
 \noindent \rem Such functions can be viewed as functions which are holomorphic in $\mathcal{W}_\delta$ with respect to the set of variables $Y$, holomorphic in a strip $\{ X\in \mathbb{C}^d, \forall j, 1\leq j\leq d, \lvert Im X_j\rvert \leq r\}$ with respect to the set of variables $X$, and $2\pi$-periodic with respect to each variable $X_j$.
 
 \bigskip
 Also, notice that for all $f\in C^\omega_{r,\delta}$, 
 
 \begin{equation}\label{supnorm}\sup_{(X,Y)\in\mathcal{V}_r\times \mathcal{W}_\delta} \lvert f(X,Y)\rvert \leq \lvert f\rvert_{r,\delta}
 \end{equation}
 
 \noindent and a classical estimate on the Fourier coefficients and a Cauchy estimate imply that for all $P\in\mathbb{Z}^d,Q\in\mathbb{N}^n$,

 \begin{equation}\label{coeffPQ}\lvert f_{P,Q}\rvert \leq \sup_{(X,Y)\in\mathcal{V}_r\times \mathcal{W}_\delta} \lvert f(X,Y)\rvert e^{-\lvert P\rvert r}\delta^{-\lvert Q\rvert }
 \end{equation}
 
 \bigskip
 \noindent 
 The space of vector fields defined on $\mathcal{V}_r\times \mathcal{W}_\delta$ with components in $C^\omega_{r,\delta}$ will be denoted by $VF_{r,\delta}$. For 
 
 $$F(X,Y)=\sum_{j=1}^d F_j(X,Y)\frac{\partial }{\partial X_j}+ \sum_{j'=1}^n F'_{j'} (X,Y)\frac{\partial }{\partial Y_{j'}} \in VF_{r,\delta}$$
 
 \noindent the norm of $F$ is by definition
 
 $$\lvert \lvert F\rvert\rvert_{r,\delta} = \max \{ \lvert F_j\rvert_{r,\delta}, \lvert F'_{j'}\rvert_{r,\delta}, 1\leq j\leq d, 1\leq j'\leq n\}$$
 
%
%
 \noindent
 For all $r,r',\delta,\delta'>0$, when dealing with an operator sending $VF_{r,\delta}$ into $VF_{r',\delta'}$, we will denote its operator norm by $\lvert \lvert \lvert \cdot \rvert\rvert\rvert_{VF_{r,\delta}\rightarrow VF_{r',\delta'}}$. 
 
 \subsection{Formal aspects}
 
 \begin{definition}
 \begin{itemize}
 \item
A formal Fourier-Taylor series in $(X,Y)$ is a series $f$ of the form

\begin{equation} \label{formalseries}f(X,Y)=\sum_{P\in\mathbb{Z}^d, Q\in \mathbb{N}^n} f_{P,Q} e^{i\langle P,X\rangle } Y^Q
\end{equation}

\noindent with coefficients $f_{P,Q}\in \mathbb{C}$ and such that for all $Q\in \mathbb{N}^n$, 

\begin{equation}\label{squaresum}\sum_{P\in \mathbb{Z}^d} \lvert f_{P,Q}\rvert ^2 <\infty
\end{equation}

\noindent i.e for all $Y\in \mathbb{C}^n$, $f(\cdot ,Y)\in L^2(\mathbb{T}^d)$
(this restriction is motivated by the fact that we need products and finite sums of Fourier-Taylor series to be themselves Fourier-Taylor series).
 
 \item
 A formal vector field (on $\mathbb{T}^d\times \mathbb{C}^n$ ) is a $d+n$-uple whose components are formal Fourier-Taylor series. 
 
\end{itemize}

 \end{definition}
 
 \begin{definition}\label{truncvf}(truncation of a formal series, of a vector field) Let $f$ be a formal series. We will denote by $T^k f$ its truncation at order $k\in \mathbb{N}$ (or $k$-jet): 

$$T^k f (X,Y)= \sum_{P\in\mathbb{Z}^d, \lvert Q\rvert \leq k}
f_{P,Q}e^{i\langle P,X\rangle } Y^Q$$

\noindent
For a vector field (resp. diffeomorphism) $F$ with components $(F_1,\dots, F_{d+n})$, its truncation at order $k$ is the vector field (resp. diffeomorphism) $T^kF$ with components

$$(T^k F_1, \dots, T^k F_d,T^{k+1}F_{d+1},\dots T^{k+1} F_{d+n})$$

\end{definition}

\noindent
\rem Thus for a vector field $F$, its truncation $T^k F$ has degree $k$ (as a vector field).

 \begin{definition}
 

 A formal diffeomorphism (of $\mathbb{C}^d/\mathbb{Z}^d\times \mathbb{C}^n$) is a sequence of analytic diffeomorphisms $(\Phi_k)$ defined on a non-increasing sequence of domains $\mathcal{V}_k\subset \mathbb{C}^d/\mathbb{Z}^d\times \mathbb{C}^n$ (for the order induced by the inclusion), such that $\Phi_k$ has degree $k$ and for all $j\leq k$, $T^j\Phi_k=\Phi_{j\lvert \mathcal{V}_k}$. We will denote by $T^k\Phi$ the analytic diffeomorphism $\Phi_k$.

   \end{definition}
 
 \bigskip
 \noindent
 \rem An analytic diffeomorphism is a particular case of a formal diffeomorphism.

 \begin{definition}(order of a formal series) A formal series $f(X,Y)= \sum_{P\in\mathbb{Z}^d, Q\in \mathbb{N}^n} f_{P,Q} e^{i\langle P,X\rangle } Y^Q$ has order $k\in \mathbb{N}$ if $\lvert Q\rvert \leq k-1\Rightarrow 
 f_{P,Q}=0$. 
 \end{definition}
 
 \noindent
 \rem Only the order in $Y$ is considered. Thus, the product of two formal series has order at least the sum of their orders. Notice that the order is not uniquely defined by maximality: if $f$ has order $k$ and $k'\leq k$, then $f$ has order $k'$.
 
 \bigskip
 \begin{definition} (order of a vector field)
For 

$$F=\sum_{j=1}^d F_j(X,Y)\frac{\partial }{\partial X_j}+ \sum_{j'=1}^n F'_{j'} (X,Y)\frac{\partial }{\partial Y_{j'}} \in VF_{r,\delta},$$

\noindent  the vector field $F$ has quasi-order $k\in \mathbb{N}$ if all the $F_j$ have order $k$ in $Y$ and the 
$F'_{j'}$ have order $k+1$ in $Y$ (as formal series).
\end{definition}

\begin{definition}
(degree of a vector field) A vector field $F$ has degree $k\in\mathbb{N}$ if its components $F_j$ have degree $k$ for $1\leq j\leq d$ and have degree $k+1$ for $d+1\leq j\leq d+n$.
\end{definition}

\noindent
\rem In particular, if a vector field has degree $k$ and order $k$, it is quasi-homogeneous of order $k$. For instance, the quasilinear vector field $S$ defined in \eqref{S} is quasi-homogeneous of order $0$.

\begin{definition}
For a vector field $F$, its quasilinear part is the vector field $T^0 F$.
\end{definition}



\begin{definition}
A formal diffeomorphism $\Phi=(\Phi_1,\dots, \Phi_{n+d})$ is tangent to identity if for all $1\leq L\leq d$, $\Phi_L(X,Y)-X_L$ has order 1 and if for all $d+1\leq L\leq d+n$, $\Phi_L(X,Y)-Y_{L-d}$ has order 2.
\end{definition}

\subsection{Resonances}

The quasilinear part $S$ defines an equivalence relation $\sim $ on $\mathbb{Z}^d\times \mathbb{N}^n$ as follows: for all $(P,Q), (P',Q')\in \mathbb{Z}^d\times \mathbb{N}^n$, 

\begin{equation}\label{equivalencerel} (P,Q)\sim (P',Q') \Leftrightarrow i\langle P,\omega\rangle+\langle Q,\Lambda\rangle =i\langle P',\omega\rangle+\langle Q',\Lambda\rangle 
\end{equation}

\begin{definition} The vectors $\omega$ and $\Lambda$ are jointly non resonant if the equivalence classes of the relation $\sim$ are reduced to singletons.
\end{definition}

\noindent
\rem It does not hold in general. Moreover, if an equivalence class is not a singleton, then it has infinitely many elements. 
In this article, $\omega$ and $\Lambda$ are \textbf{not assumed to be jointly non resonant.}

\bigskip
\noindent 
For $c\in\mathbb{C}$, $\mathcal{C}_c$ denotes the equivalence class such that for all $(P,Q)\in\mathcal{C}_c$, $i\langle P,\omega \rangle +\langle Q,\Lambda\rangle =c$.

\bigskip
\noindent Let $f$ be an analytic function. Then the Fourier-Taylor development of $f$ has only indices in $\mathcal{C}_c$ if and only if $S(f)=cf$, where $S(f)$ is the Lie derivative of $f$ along $S$:

$$S(f)= \sum_{j=1}^d \omega_j \frac{d}{dX_j} f + \sum_{j'=1}^n \lambda_{j'} Y_{j'} \frac{d}{dY_{j'}} f$$

\begin{definition}
\begin{itemize}
\item  Monomials $f_{P,Q}e^{i\langle P,X\rangle } Y^Q$ with $i\langle P,\omega\rangle +\langle Q,\Lambda\rangle =0$ are resonant, or invariant, monomials (indeed they are constant in the direction of the vector field $S$).
\item 
For a formal series $f$, its resonant part is the sum of all its resonant monomials, and its non resonant part is the sum of all other monomials. 

\item
A formal series is resonant if it is equal to its resonant part, and non resonant if it is equal to its non resonant part.
\end{itemize}\end{definition}

\noindent
\rem This means that a formal series is resonant if all its monomials are indexed by $(P,Q)\in\mathcal{C}_0$, i.e if $S(f)=0$.

\begin{definition}

Let $F=\sum_{j=1}^d F_j\frac{\partial }{\partial X_j}+\sum_{j'=1}^nF_{d+j'}\frac{\partial }{\partial Y_{j'}}$ be a vector field on a domain in $\mathbb{C}^d/\mathbb{Z}^d\times \mathbb{C}^n$. 
\begin{itemize}
\item $F$ is non resonant if $F_1,\dots, F_d$ are non resonant formal series and if for all $1\leq j\leq n$, $F_{d+j}(X,Y) Y_j^{-1}$ is a non resonant formal series.


\item $F$ is resonant if $F_1,\dots, F_d$ are resonant formal series and if for all $1\leq j\leq n$, $F_{d+j}(X,Y) Y_j^{-1}$ is a resonant formal series.
\item There is a unique decomposition $F=F_{res}+F_{nr}$ where $F_{res}$ is a resonant vector field and $F_{nr}$ is a non resonant vector field; then $F_{res}$ is called the resonant part of $F$ and $F_{nr}$ its non resonant part.
\end{itemize}

\end{definition}

\noindent
\rem Let $j$ with $1\leq j\leq n$ and write the decomposition 

$$F_{d+j}(X,Y)= \sum_{P,Q} F_{d+j,P,Q} e^{i\langle P,X\rangle} Y^Q$$

\noindent Then $Y_j^{-1} F_{d+j}(X,Y)$ is resonant if 

$$F_{d+j,P,Q}\neq 0\Rightarrow i\langle P,\omega\rangle +\langle Q,\Lambda\rangle -\lambda_j=0$$ 

\noindent
This coincides with Aurouet's definition of the invariants.

\begin{definition} \label{resonantdiffeo}
Let $\Phi$ be a formal diffeomorphism. Suppose that for all $1\leq j\leq d$, $\Phi_j(X,Y)=X_j+\tilde{\Phi}_j(X,Y)$ where $\tilde{\Phi}_j(X,Y)$ has order 1, and that for all $1\leq j'\leq n$, $\Phi_{d+j'}(X,Y)=Y_{j'} \tilde{\Phi}_{d+j'}(X,Y)$ where $ \tilde{\Phi}_{d+j'}(X,Y)$ has order 1.
\begin{itemize}
\item $\Phi$ is non resonant up to order $k$ if for all $1\leq j\leq n+d$, $T^k\tilde{\Phi}_j$ is non resonant;

\item  $\Phi$ is resonant up to order $k$ if for all $1\leq j\leq n+d$, $T^k\tilde{\Phi}_j$ is resonant;
\item $\Phi$ is non resonant (resp. resonant) if for all $k\in \mathbb{N}$ it is non resonant up to order $k$ (resp. resonant up to order $k$).
\end{itemize}
\end{definition}

\bigskip 
\noindent
\rem This definition applies in particular if $\Phi$ is an analytic diffeomorphism.

\begin{definition}
Let $F$ be a formal vector field and $G$ be an analytic vector field; $G$ is formally conjugate to $F$ if there exists a formal diffeomorphism $\Phi$ such that for all $k\in\mathbb{N}$, 

$$T^k(DT^k\Phi\cdot F)=T^k(G\circ T^k\Phi)$$

\end{definition}

\noindent
\rem If an analytic vector field $G$ is formally conjugate to a formal vector field $F$ and if $G$ is also analytically conjugate to an analytic vector field $H$, then $H$ is formally conjugate to $F$. Indeed,
assume that for all $k\in\mathbb{N}$

$$T^k(DT^k\Phi\cdot F)=T^k(G\circ T^k\Phi)$$

\noindent and that there exists an analytic diffeomorphism $\tilde{\Phi}$ such that 

$$D\tilde{\Phi}\cdot G=H\circ \tilde{\Phi}.$$

\noindent Let $\Psi_k=T^k(\tilde{\Phi}\circ T^k\Phi)$; 
then $D\Psi_k = T^k( D  \tilde{\Phi} \circ T^k\Phi \cdot D T^k\Phi)$;
thus for all $k$,

\begin{equation}\begin{split}T^k(D\Psi_k\cdot F)&= T^k [D  \tilde{\Phi} \circ T^k\Phi \cdot D T^k\Phi \cdot F]=T^k [D \tilde{\Phi}\circ T^k\Phi \cdot G\circ T^k{\Phi}]\\
&=T^k[H\circ \tilde{\Phi}\circ T^k\Phi]=T^k[H\circ \Psi_k]
\end{split}\end{equation}

\noindent therefore the sequence $(\Psi_k)$ defines a formal diffeomorphism conjugating $H$ to $F$.

\bigskip
\noindent \textbf{Other notations:} For $I=(i_1,\dots, i_{d+n})\in \mathbb{N}^{d+n}$, the symbol $\partial_I$ will be a compact notation for the partial derivative $\partial_{X_1}^{i_1}\dots \partial_{X_d}^{i_d} \partial_{Y_1}^{i_{d+1}}\dots \partial_{Y_n}^{i_{d+n}}$.

\noindent
Sometimes the remainder of order $k$ in $Y$ in a formal series or a formal vector field will be denoted by $O(Y^k)$.

\noindent
If $F\in VF$ and $g$ is a function, the Lie derivative of $g$ with respect to $F$ will be denoted by $F(g)$.

\section{Setting and arithmetical conditions}\label{setup}

Let $r_0>0, \delta_0>0$.
Consider the following system on the neighborhood $\mathcal{V}_{r_0}\times \mathcal{W}_{\delta_0}\subset \mathbb{C}^d/ \mathbb{Z}^d \times \mathbb{C}^n$ of $\mathbb{T}^d\times \{0\}$:

\begin{equation}\label{chres}
\left\{ \begin{array}{l}
\dot{X}_1=\omega_1 +R_1(X,Y)\\
\vdots\\
\dot{X}_d=\omega_d+R_d(X,Y)\\
\dot{Y}_1=\lambda_1 Y_1 +R_{d+1}(X,Y)\\
\vdots\\
\dot{Y}_n=\lambda_n Y_n +R_{d+n}(X,Y)\\
\end{array}\right.
\end{equation}

\noindent 
where $(X_1,\dots, X_d)\in \mathcal{V}_{r_0}$, $(Y_1,\dots ,Y_n)\in  \mathcal{W}_{\delta_0}$, $R_1,\dots, R_d$ have order 1 and $R_{d+1},\dots, R_{d+n}$ have order 2. Suppose that 
$R_1,\dots, R_{d+n}$ are all in $C^\omega_{r_0,\delta_0}$.

\noindent
Denote by $S$ the quasilinear part of the vector field generating the system \eqref{chres}:

$$S(Y)=\sum_{j=1}^d \omega_j\frac{\partial }{\partial X_j}+\sum_{j'=1}^n \lambda_{j'} Y_{j'} \frac{\partial }{\partial Y_{j'}}$$

\noindent
Thus for all $r,\delta >0$, $\lvert \lvert S\rvert\rvert_{r,\delta} =\max \{ \omega_j, \lambda_{j'} \delta, 1\leq j\leq d,1\leq j'\leq n\}$ and if $\delta\leq \min_{j'}\max_j\lvert \omega_j / \lambda_{j'}\rvert$, then $\lvert \lvert S\rvert\rvert_{r,\delta}\leq \max_j \lvert \omega_j\rvert:=C_\omega$. For short, denote $F=S+R$ where $R$ is the vector field with components $R_1,\dots, R_{d+n}$ defined in \eqref{chres}.

%

\noindent
The following arithmetical condition, which will be assumed throughout the article, is a reformulation of Brjuno's $\gamma$ and $\omega$ conditions (also used in \cite{Au13}). It will be used in the theorem of analytic normalization.

\begin{ass}\label{gamma+omega} There exists a positive increasing unbounded function $g$ on $[1,+\infty[$, 
an increasing sequence of integers $(m_k)_{k\geq 0}$ and two sequences of positive real numbers $(\epsilon_k)_{k\geq 0},(r_k)_{k\geq 0}$ such that:
\begin{enumerate}
\item \label{mkcroiss} $m_0=1;m_{k+1}\leq 2m_k+1$;
\item \label{brunosum} $B:=\sum_{k\geq 1} \frac{\ln g(m_k)}{m_k}<+\infty$,
\item \label{convergence} $\forall k\geq 0, (n+2)g(m_k)\ln g(m_k) \geq m_k$,
\item \label{highfq} $\forall P\in \mathbb{Z}^d, \forall Q\in \mathbb{N}^n$, 

$\lvert P\rvert > m_k,\lvert Q\rvert \leq m_k, i\langle P,\omega\rangle +\langle Q,\Lambda\rangle \neq 0  \Rightarrow  \lvert i\langle P,\omega\rangle +\langle Q,\Lambda\rangle \rvert^{-1} \leq e^{\epsilon_k \lvert P\rvert}$,
\item \label{lowfq} $\forall P\in \mathbb{Z}^d,\forall Q\in \mathbb{N}^n$, 

$\lvert P\rvert \leq m_k, \lvert Q\rvert \leq m_k,i\langle P,\omega\rangle +\langle Q,\Lambda\rangle \neq 0 \Rightarrow \lvert i\langle P,\omega\rangle +\langle Q,\Lambda\rangle \rvert^{-1} \leq g(m_k)$,
\item \label{rkassezgrand} $\forall k, 4^{4d+6}(r_k-r_{k+1}-\epsilon_k)^{-4d-6} \leq g(m_{k})$,
\item \label{poslimitrk} $\forall k\geq 0, r_k >\frac{1}{2}$.
\end{enumerate}
\end{ass}

\noindent 
\rem
The case when $m_k=2^k$ corresponds to Brjuno's standard $\omega$ condition. However, the assumption \ref{gamma+omega}.\ref{mkcroiss} seems to prevent $(m_k)$ from increasing faster than $2^k$.

\bigskip
\noindent
Note that the set of elements $g,(m_k),(\epsilon_k),(r_k)$ satisfying assumption \ref{gamma+omega} is non-empty. In \cite{Au13}, the conditions $\gamma$ and $\omega$ correspond to the case $m_k=2^k, r_k=r_0\prod_{j=1}^k \frac{1}{g(m_j)^{\frac{1}{2^j}}2^{\frac{2j}{2^j}}}, \epsilon_k>\frac{1}{4} (1-\frac{1}{2^{\frac{k}{2^k}}})$, and $g(m_k)=\max \{ (i\langle P,\omega\rangle +\langle Q,\Lambda\rangle )^{-1} / \lvert P\rvert \leq m_k, \lvert Q\rvert \leq m_k\}$, with $r_0$ big enough in order to have $r_k>\frac{1}{2}$.

\bigskip
\noindent The case where $g(m_k)=C m_k^\tau$ for $m_k=2^k$, for a positive constant $C$ and an exponent $\tau\geq d+n+1$, corresponds to the situation where the vector $(\omega,\Lambda)$ satisfies a Siegel diophantine condition; therefore the set of $(\omega,\Lambda)\in \mathbb{C}^{d+n}$ satisfying the assumption \ref{gamma+omega} has a large Lebesgue measure.

\bigskip
\noindent Also, assumption \ref{gamma+omega} does not contain any upper bound on $g(m_0)$: the only upper bound on $g$ comes from \ref{brunosum} and it is logarithmic; therefore multiplying $g(m_k)$ by a constant does not change the condition. 
In section \ref{mainthm}, we shall assume that for all $k\geq 0$, $g(m_k)$ is larger than a fixed constant.


%
%
%
%
%
%
%

\section{Formal normal form}


This section contains an explicit construction of a formal normal form for the system \eqref{ptdepart}. A direct method will be used, i.e the formal conjugation will be written explicitly, by a Taylor development, as a recurrence relation between two consecutive orders of truncation of the change of variables. The following lemma explains why the recurrence can be solved.

\begin{lem}\label{opdiff}Let $\bar{\Phi}_1,\dots ,\bar{\Phi}_d$ be analytic Fourier-Taylor series in $(X,Y)\in\mathbb{T}^d\times \mathbb{C}^n$ which are of order 1 in $Y$ and $\bar{\Phi}_{d+1},\dots ,\bar{\Phi}_{d+n}$ be analytic Fourier-Taylor series in $(X,Y)\in\mathbb{T}^d\times \mathbb{C}^n$ which are of order 2 in $Y$. 
Let $I=(i_1,\dots, i_{n+d})\in\mathbb{N}^{d+n}$. For any analytic vector field $F$, let

\begin{equation}D_IF(X,Y)=c_I
\partial_IF(X,Y)
\cdot \bar{\Phi}_1(X,Y)^{i_1}\dots \bar{\Phi}_{d+n}(X,Y)^{i_{d+n}} 
\end{equation}

\noindent where $c_I\in \mathbb{C}^*$ and $\partial_I$ is a compact notation for $\partial_{X_1}^{i_1}\dots \partial_{X_d}^{i_d}\partial_{Y_1}^{i_{d+1}}\dots \partial_{Y_n}^{i_{d+n}}$. If $F$ has order $k$ in $Y$, then $D_IF$ has order $k+i_1+\dots +i_{d+n}$.

%

\end{lem}

\dem It is enough to prove the statement for elementary $(i_1,\dots, i_{d+n})$. For $1\leq j\leq d$, the derivative $\partial_{X_j} $ does not change the order in 
$Y$ and multiplying by $\bar{\Phi}_j$ adds a degree in $Y$. Thus the order of $D_{e_j}F$ is $k+1$ if $F$ has order $k$.

%
%
%
%
%

For $1\leq j\leq n$ and $1\leq l\leq d$,
the derivative $\partial_{Y_j}F_l$ has order $k-1$ and multiplying by $\bar{\Phi}_{d+j}$ adds two degrees in $Y$, thus $D_{e_{d+j}}F$ has order $k+1$; 
for $d+1\leq l\leq d+n$, $F_l$ has order $k+1$, therefore the derivative $\partial_{Y_j}F_l$ has order $k$, and multiplying by $\bar{\Phi}_{d+j}$ adds two degrees in $Y$. 

\noindent Finally, $D_{e_{d+j}}F$ has order $k+1$ (as a vector field).
$\Box$

%
%
%
%
%


\bigskip
\noindent The following Proposition gives a construction of a formal normal form:

\begin{prop}
The system \eqref{chres} can be formally normalized, i.e there exist a formal diffeomorphism $\Phi$ and a formal resonant vector field $NF$ such that for all $k$,

$$T^k[ D(T^k\Phi) \cdot NF]=T^k(F\circ T^k \Phi)=T^k((S+R)\circ T^k\Phi)$$

\end{prop}

\dem 
One looks for a formal diffeomorphism $\Phi$ tangent to identity, i.e such that 

$$T^0(\Phi-Id)=0$$

\noindent and for a resonant vector field $NF=S+N$ where $N$ has order 1, such that
%
%
%
%

\begin{equation}\label{conjordrek}T^k (DT^k\Phi\cdot NF) = S\circ T^k \Phi+ T^k (R\circ T^k\Phi) 
\end{equation}

%
%
%
%

\noindent
The formal diffeomorphism $\Phi$ and the normal form $NF$ will be constructed gradually, every quasi-homogeneous part of degree $k$ being given formally by the data $F$ and by the truncations $T^{k-1}NF$ and $T^{k-1}\Phi$. 

\bigskip
\noindent Developing the composition $R\circ T^k \Phi$
in Taylor series 
(recall that $F$ is analytic),
letting $\bar{\Phi}=\Phi-Id$
and developing in $\bar{\Phi}$, 
one obtains

%
\begin{equation}\begin{split}
R\circ (Id +T^k\bar{\Phi})&=R+\sum_{I=(i_1,\dots, i_{d+n})\neq (0,\dots, 0)} c_I 
\partial_I R\cdot 
T^k\bar{\Phi}_1^{i_1}\dots T^k\bar{\Phi}_{d+n}^{i_{d+n}}
\end{split}
\end{equation}

%

\noindent where $c_I=\frac{1}{i_1!\dots i_{d+n}!}$ and $\partial_I$ is the partial derivative $\partial_{X_1}^{i_1}\dots \partial_{Y_n}^{i_{d+n}}$, 
that is to say,

\begin{equation}\label{blabla}R\circ (Id +T^k\bar{\Phi})(X,Y)=R(X,Y)+\sum_{(i_1,\dots, i_{d+n})\neq (0,\dots ,0)} D_{i_1,\dots, i_{n+d}}R(X,Y)
\end{equation}

\noindent where $D_{i_1,\dots, i_{n+d}}$ is as in Lemma \ref{opdiff} as long as the $\bar{\Phi}_j$ have order 1 for $j\leq d$ and order 2 for $j\geq d+1$, i.e as long as $\Phi$ is tangent to identity; the operator $D_{i_1,\dots, i_{n+d}}$ therefore increases the order by $i_1+\dots +i_{n+d}$. Although the sum in \eqref{blabla} does not necessarily define a formal series since it has an infinite number of terms, its truncation has a finite number of terms, therefore it defines a formal series.

Thus if $\Phi$ is tangent to identity, the conjugation \eqref{conjordrek} can be rewritten

\begin{equation}\label{TkLleqd}\begin{split}&T^k NF + T^k(DT^k\bar{\Phi}\cdot NF)\\
& = S\circ T^k\Phi+T^kR + \sum_{(i_1,\dots, i_{d+n})\neq (0,\dots, 0)} T^k (D_{i_1,\dots, i_{n+d}}T^{k-i_1\dots -i_{n+d}}R)
\end{split}\end{equation}

\bigskip
\noindent
Let $k\geq 1$
; suppose that $T^{k-1}NF$ and $T^{k-1}{\Phi}$ are known and analytic, and that $T^0\bar{\Phi}=0$. 
In Equation \eqref{TkLleqd}, 
%
%
%
the sum in the right-hand side is given because $R$ has order 1 and therefore for all $I$, 

$$D_I T^{k-\lvert I\rvert}R=c_I\partial_I T^{k-\lvert I\rvert } R\cdot T^{k-1}\bar{\Phi}^I + O(Y^k)$$

On the left-hand side, since the $N_j$ have order 1 for $1\leq j\leq d$ and order 2 for $d+1\leq j\leq d+n$, then 

$$T^k(DT^k\bar{\Phi}\cdot N)
= T^k\left(\sum_j\partial_{X_j}T^{k-1}\bar{\Phi}N_j+\sum_{j'}\partial_{Y_{j'}} T^{k-1}\bar{\Phi} N_{j'+d}\right)$$

\noindent
and since $T^0\bar{\Phi}=0$, then 

\begin{equation}\label{quantity}T^k(DT^k\bar{\Phi}\cdot N)
= T^k\left(\sum_j\partial_{X_j}T^{k-1}\bar{\Phi}T^{k-1}N_j+\sum_{j'}\partial_{Y_{j'}} T^{k-1}\bar{\Phi} T^kN_{j'+d}\right)
\end{equation}

Note that $T^{k-1}N_j$ and $T^kN_{j'+d}$ are exactly the components of $T^{k-1}N$, which is assumed to be a known analytic vector field.
Since $T^{k-1}\bar{\Phi}$ is also known and analytic, 
then the whole quantity \eqref{quantity} is given and analytic.
Thus, in \eqref{TkLleqd}, only $T^k N-T^{k-1}N + T^k(D\bar{\Phi}\cdot S)-T^{k-1}(D\bar{\Phi}\cdot S)-(S\circ T^k\Phi -S\circ T^{k-1}\Phi)$ is unknown. One therefore has to solve

$$T^k N-T^{k-1}N + T^k(D\bar{\Phi}\cdot S)-T^{k-1}(D\bar{\Phi}\cdot S)-(S\circ T^k\Phi -S\circ T^{k-1}\Phi)=G$$

where $G$ is a given analytic vector field of degree $k$, namely 

$$G=\sum_I T^k(D_I T^{k-\lvert I\rvert } R)-S-T^{k-1}N-T^k(DT^k\bar{\Phi}\cdot N) -T^{k-1}D\bar{\Phi}\cdot S+S\circ T^{k-1}\bar{\Phi}$$

\noindent
which is homogeneous in $Y\in \mathcal{W}_{\delta_0}$, of degree $k$, by construction of $T^{k-1}N$ and $T^{k-1}\Phi$. The Fourier-Taylor coefficient $G_{P,Q}$ is well defined since the sum runs over a finite number of indices.


\noindent Letting 

$$T^k N-T^{k-1}N=G_{res}$$

\noindent and

$$\bar{\Phi}_{L,P,Q}=(i\langle P,\omega\rangle +\langle Q,\Lambda\rangle)^{-1}G_{L,P,Q}$$

\noindent for non resonant $P,Q $ such that $\lvert Q\rvert =k$ and $L\leq d$, and 

$$\bar{\Phi}_{L,P,Q}= (i\langle P,\omega\rangle +\langle Q,\Lambda\rangle -\lambda_{L-d} )^{-1}G_{L,P,Q}$$

\noindent for all $P,Q$ such that $i\langle P,\omega\rangle +\langle Q,\Lambda\rangle -\lambda_{L-d}\neq 0$, with $\lvert Q\rvert =k$ and $L\geq d+1$, the formal diffeomorphism $\Phi$ and of the normal form $NF$ are defined up to order $k$. $\Box$

\section{Analytic normal form}

This section deals with the problem of the existence of an analytic normalization.

%

%
%
%
%

\subsection{Preliminary observations}

The following is a lemma stating that being resonant is preserved by composition with a resonant diffeomorphism.

\begin{lem}\label{resonantcomp}
Let $\Phi$ be a resonant holomorphic diffeomorphism from $\mathcal{V}\times \mathcal{W}\subset \mathbb{C}^d/\mathbb{Z}^d\times \mathbb{C}^n$ to $\mathcal{V}'\times \mathcal{W}'\subset \mathbb{C}^d/\mathbb{Z}^d\times \mathbb{C}^n$, which is tangent to identity.
%
%
%

\noindent
 If $R $ is an analytic resonant vector field on $\mathcal{V}'\times \mathcal{W}'$
 , then $R\circ \Phi$ is also resonant. 
\end{lem}

\dem Developing in Taylor series, for all $(X,Y)\in \mathcal{V}\times \mathcal{W}$,

$$R\circ\Phi(X,Y)= \sum_{I\in\mathbb{N}^{n+d}} c_I
\partial_I R(X,Y)
[(\Phi-Id)(X,Y)]^I
$$

\noindent with $c_I =\frac{1}{I_1!\dots I_{d+n}!}$.
Now $\Phi$ is resonant, which implies (see Definition \ref{resonantdiffeo}) that 

$$(\Phi-Id)(X,Y)^I Y_1^{-i_{d+1}}\dots Y_n^{-i_{d+n}}$$

\noindent is a resonant series. 



Thus, for all $I$ in the sum, $\partial_I R(X,Y) [(\Phi-Id)(X,Y)]^I$ is a resonant vector field, and this property passes to the sum. $\Box$

%
%
%
%

\begin{lem}\label{prodres}
Let $R$ be a formal resonant vector field and $\Phi$ a resonant analytic diffeomorphism which is tangent to identity. Then $D\Phi\cdot R$ is a formal resonant vector field. 
\end{lem}

\dem 
Let $\tilde{\Phi}_1,\dots, \tilde{\Phi}_{d+n}$ be such that 

$$\Phi=(X_1+\tilde{\Phi}_1,\dots, X_d+\tilde{\Phi}_d, Y_1\tilde{\Phi}_{d+1},\dots,Y_n\tilde{\Phi}_{d+n})$$

\noindent then 

\begin{equation}\begin{split}
&D\Phi\cdot R=\\
&\left(\begin{array}{c}
R_1+\partial_{X_1}\tilde{\Phi}_1R_1+\dots +\partial_{X_d}\tilde{\Phi}_1R_d+ \partial_{Y_1}\tilde{\Phi}_1R_{d+1}+\dots + \partial_{Y_n}\tilde{\Phi}_1R_{d+n}\\
\vdots\\
R_d+\partial_{X_1}\tilde{\Phi}_dR_1+\dots +\partial_{X_d}\tilde{\Phi}_dR_d+ \partial_{Y_1}\tilde{\Phi}_dR_{d+1}+\dots + \partial_{Y_n}\tilde{\Phi}_dR_{d+n}\\
Y_1\sum_{j=1}^d\partial_{X_j}\tilde{\Phi}_{d+1}R_j
+ \tilde{\Phi}_{d+1} R_{d+1}+ Y_1\sum_{j'=1}^n\partial_{Y_{j'}} \tilde{\Phi}_{d+1}R_{d+j'}
\\
\vdots\\
Y_n\sum_{j=1}^d\partial_{X_j}\tilde{\Phi}_{d+n}R_j
+  \tilde{\Phi}_{d+n}R_{d+n} +Y_n\sum_{j'=1}^n\partial_{Y_{j'}}\tilde{\Phi}_{d+n}R_{d+j'}\\
\end{array}\right)
\end{split}\end{equation}

\noindent Since all the $\tilde{\Phi}_j$ are resonant by assumption, and since $R$ is a resonant vector field, then $D\Phi\cdot R$ is a resonant vector field. $\Box$

\bigskip
\noindent
The following proposition states that if a system has been normalized up to some order $k$, then the formal normalizing diffeomorphism of this partially normalized system has to be resonant up to order $k$.

\begin{prop}\label{resonantorderk}
Let $\Phi$ be a formal diffeomorphism.
%
%
%

\noindent
Suppose that $\Phi$ formally conjugates 
an analytic vector field $S+N_k+R_k$ where $N_k $ is resonant of order 1 and $R_k$ has order $k+1$,
to a formal resonant vector field $NF_1=S+R$, where $R$ has order 1
: for all $j\in\mathbb{N}$,

$$T^j(DT^j\Phi\cdot NF_1)=T^j[(S+N_k+R_k)\circ T^j\Phi]$$

\noindent 
Then $T^k\Phi$ is a resonant diffeomorphism.
\end{prop}

\dem 
By assumption and since $R_k$ has order $k+1$, for all $l\leq k$,

$$T^l(DT^l\Phi\cdot NF_1)=S\circ T^l\Phi+T^l(N_{k}\circ T^l \Phi )$$

\noindent 
and since $N_k$ has order 1, 

$$T^l(DT^l\Phi\cdot NF_1)=S\circ T^l\Phi+T^l(N_{k}\circ T^{l -1}\Phi )$$

\noindent 
Now $NF_1=S+R$ where $R$ has order 1, thus

$$T^l(DT^l\Phi\cdot S)+T^l(DT^{l-1}\Phi\cdot R)=S\circ T^l\Phi+T^l(N_{k}\circ T^{l -1}\Phi )\cdot $$

%
%
%
%
%
%
%
%
%

\noindent
By induction, let $l\leq k-1$ and assume that $T^{l-1} {\Phi}$ is a resonant diffeomorphism. By Lemma \ref{resonantcomp}, this implies that $T^l(N_{k}\circ T^{l -1}\Phi)$ is a resonant vector field. 
Moreover, $T^{l}(DT^{l-1}{\Phi}\cdot R )$ is resonant by Lemma \ref{prodres}.

%
\noindent Thus, $T^l(DT^l{\Phi}\cdot S)-S\circ T^l\Phi$ is a resonant vector field. 

\noindent
Let $j\leq d$. Since

$$T^l(DT^l{\Phi}_j\cdot S )-\omega_j= \sum_{P,\lvert Q\rvert \leq l} T^l{\Phi}_{j,P,Q} (i\langle P,\omega\rangle +\langle Q,\Lambda\rangle ) e^{i\langle P,X\rangle } Y^Q-\omega_j$$

\noindent
and since this quantity is a resonant function as shown above,
then for all non-resonant $(P,Q)$ with $\lvert Q\rvert \leq l$, ${\Phi}_{j,P,Q}=0$ and therefore ${\Phi}_j$ is resonant up to order $l$. 
Similarly, if $j\geq d+1$, 

$$T^l(DT^l{\Phi}_j\cdot S )-\lambda_{j-d}T^l\Phi_j= \sum_{P,\lvert Q\rvert \leq l} {\Phi}_{j,P,Q} (i\langle P,\omega\rangle +\langle Q,\Lambda\rangle-\lambda_{j-d} ) e^{i\langle P,X\rangle } Y^Q$$

\noindent
and since $Y_{j-d}^{-1}(T^l(DT^l{\Phi}_j\cdot S )-\lambda_{j-d}T^l\Phi_j)$ is resonant, then ${\Phi}_{j,P,Q}=0$ if $i\langle P,\omega\rangle +\langle Q,\Lambda\rangle-\lambda_{j-d}\neq 0$ for all $\lvert Q\rvert \leq l$, which implies that $T^l\Phi_j(X,Y)Y_{j-d}^{-1}$ is a resonant function.

\bigskip

\noindent Therefore, $T^l\Phi$ is a resonant diffeomorphism.
By immediate induction, the formal diffeomorphism ${\Phi}$ is resonant up to order $k$. $\Box$

%
%
%
%
%
%
%
%
%

\bigskip
The next proposition states that the algebraic property \eqref{hypsurFN} is preserved by partial normalization tangent to identity; in other terms, if a vector field has been normalized analytically up to order $k$, i.e if it is analytically conjugate to a vector field $N_k+R_k$ where $N_k$ is resonant and $R_k$ has order $k+1$, and if the same vector field is formally conjugate to a resonant vector field of the form \eqref{hypsurFN}, then the partial normal form $N_k$ also satisfies \eqref{hypsurFN}.

\begin{prop}\label{nfk=aS}
Let $F_2$ be an analytic vector field on $\mathcal{V}\times \mathcal{W}$ with quasilinear part $S$. Assume that $F_2$ is formally normalizable, i.e there exists a formal diffeomorphism $\Phi$ which is tangent to identity and a formal vector field $N_1$ such that $\Phi$ conjugates $F_2$ to $N_1$: for all $k\in\mathbb{N}$,

\begin{equation}\label{normpartielle}T^k(DT^k\Phi\cdot N_1)=T^k(F_2\circ T^k \Phi)
\end{equation}

\noindent
 For all $k\in\mathbb{N}$, if $T^kF_2$ is a resonant analytic vector field and if 
 $T^kN_1=a_k S$ where $a_k$ is $L^2$ in the first set of variables $X$ and has degree $k$ in $Y$, with constant part 1, then $T^k F_2= bS$, for a resonant analytic function $b$ which has constant part 1.
\end{prop}

\dem 
 By \eqref{normpartielle} and by assumption on $N_1$, one has

\begin{equation}\label{conjnf}T^k(DT^k\Phi\cdot a_kS)= T^k(F_{2}\circ T^k\Phi)
\end{equation}

\noindent
Now let $\bar{\Phi}=\Phi-Id$. If $1\leq l\leq d$, 

\begin{equation}\label{X}\begin{split}
&DT^k\Phi_l(X,Y)\cdot a_k(X,Y)S\\
&=a_k(X,Y) DT^k\Phi_l(X,Y)\cdot S\\
&=a_k(X,Y)[ \sum_{j=1}^d \partial_{X_j}T^k\Phi_l(X,Y) \cdot \omega_j+ \sum_{j'=1}^n \partial_{Y_{j'}} T^k\Phi_l(X,Y)\cdot \lambda_{j'}Y_{j'}]\\
& = a_k(X,Y)[  \sum_{j=1}^d \delta_{j,l}\omega_j + \partial_{X_j}T^k\bar{\Phi}_l(X,Y)\cdot \omega_j + \sum_{j'=1}^n \partial_{Y_{j'}}T^k\bar{ \Phi}_l(X,Y)\cdot \lambda_{j'}Y_{j'}]\\
&= a_k(X,Y) [ \omega_l + \sum_{P,\lvert Q\rvert \leq k} \bar{\Phi}_{l,P,Q}e^{i\langle P,X\rangle } Y^Q(i\langle P,\omega\rangle +\langle Q,\Lambda\rangle )]
\end{split}\end{equation}

\noindent 
and the resonant part is merely $a_k(X,Y)\omega_l$. If $d+1\leq l \leq d+n$,

\begin{equation}\begin{split}\label{nfkY}
&DT^k\Phi_l(X,Y)\cdot a_k(X,Y)S\\
&=a_k(X,Y)[ \sum_{j=1}^d \partial_{X_j}T^k\Phi_l(X,Y) \cdot \omega_j+ \sum_{j'=1}^n \partial_{Y_{j'}}T^k \Phi_l(X,Y)\cdot \lambda_{j'}Y_{j'}]\\
& = a_k(X,Y) [ \sum_{j=1}^d \partial_{X_j}T^k\bar{\Phi}_l(X,Y) \cdot \omega_j+  \sum_{j'=1}^n (\partial_{Y_{j'}} T^k\bar{\Phi}_l(X,Y)+ \delta_{j',l})\lambda_{j'}Y_{j'}]\\
&=a_k(X,Y) [ \sum_{P,\lvert Q\rvert \leq k} \bar{\Phi}_{P,Q} e^{i\langle P,X\rangle} Y^Q (i\langle P,\omega\rangle +\langle Q,\Lambda\rangle ) + \lambda_lY_l]
\end{split}\end{equation}

\noindent therefore the resonant part of $Y_l^{-1} DT^k\Phi_l(X,Y)\cdot a_k(X,Y)S$ is $\lambda_l a_k(X,Y)$.
On the other side, if one assumes that $T^kF_2$ is resonant, then

$$T^k(F_{2}\circ T^k\Phi)= T^k(N_k\circ T^k\Phi +R_k\circ T^k\Phi)$$

\noindent 
where $R_k$ has order $k+1$ and $N_k$ is resonant. Since the composition increases the order, 

$$T^k(F_{2}\circ T^k\Phi)= T^k(N_k\circ T^k\Phi)$$

\noindent
By Proposition \ref{resonantorderk}, $T^k\Phi$ is resonant, therefore, by Lemma \ref{resonantcomp}, $T^k(N_k\circ T^k\Phi)$
is resonant, which implies that $T^k(F_{2}\circ T^k\Phi)$ is resonant. Thus, 
keeping only the resonant parts in \eqref{X} and \eqref{nfkY} and substituting them in the resonant part of \eqref{conjnf}, one has

%
%

$$ a_k\cdot S=T^k(N_{k}\circ T^k\Phi)
$$
%
%
%
%
%

%
%

\noindent Now this implies that $T^k N_k$ itself is proportional to $S$: indeed, letting 

$$\Phi_X(X,Y)=(T^k\Phi_1(X,Y),\dots, T^k\Phi_d(X,Y))$$ 

\noindent and 

$$\Phi_Y(X,Y)=(T^{k+1}\Phi_{d+1}(X,Y),\dots, T^{k+1}\Phi_{d+n}(X,Y))$$

\noindent (note that $\Phi_X$ and $\Phi_Y$ are analytic)
then developing in Taylor series with respect to the second set of variables $Y$, for all $0\leq j\leq k$,

\begin{equation}\begin{split}&T^j(N_k\circ (\Phi_X(X,Y),\Phi_Y(X,Y)))\\
&= \sum_{I\in\mathbb{N}^n,\lvert I\rvert\leq j} T^j \partial_IN_k(\Phi_X(X,Y),Y)\cdot (\Phi_Y(X,Y)-Y)^I
\end{split}\end{equation}

\noindent where $\lvert I\rvert =i_1+\dots +i_{n}$ for all $I=(i_1,\dots, i_{n})$ and $\partial_I=\partial_{Y_1}^{i_1}\dots \partial_{Y_n}^{i_{n}}$ (recall that by Lemma \ref{opdiff}, $f\mapsto \partial_I f \cdot (\Phi_Y(X,Y)-Y)^I$ 
increases the order by $\lvert I\rvert$). By a simple recurrence on the order of truncation $j$, one finds that 

\begin{itemize}

\item $T^0N_k=a_0\cdot S$,
\item for all $0\leq j\leq k-1$, 

\begin{equation}\begin{split}&T^jN_k(X,Y) = a_j(X,Y)\cdot S(Y) \\
&- \sum_{\lvert I\rvert =1}^j T^j(\partial_I T^{j-\lvert I\rvert}(N_k(\Phi_X(X,Y),Y)) \cdot (\Phi_Y(X,Y)-Y)^I
\end{split}\end{equation}

\noindent and if by assumption $T^{j-1}N_k =b_{j-1}S$ for an analytic function $b_{j-1}$, then for all $I$ in the sum, $T^{j-\lvert I\rvert}(N_k(\Phi_X(X,Y),Y))$ is proportional to $S$, therefore 

$${T^j(\partial_I T^{j-\lvert I\rvert}(N_k(\Phi_X(X,Y),Y)) \cdot (\Phi_Y(X,Y)-Y)^I}$$ 

\noindent
is also proportional to $S$ since 
$(\Phi_Y(X,Y)-Y)^I$ is a scalar, and it is analytic in $(X,Y)$.
\end{itemize}

\noindent 
Finally,
$T^kN_k=b\cdot S$ for an analytic function $b$. $\Box$

\subsection{Proof of Theorem \ref{mainresult}}\label{mainthm}

\noindent This section is dedicated to the proof of the main result. Let $\delta_0$ be as in Section \ref{setup}; assume that $\delta_0\leq \min(1, C_\omega)$ where $C_\omega$ was defined in Section \ref{setup}.

\begin{thm}\label{mainthmbis}
Assume that the analytic vector field $S+R\in VF_{r_0,\delta_0}$ generating \eqref{chres} is formally conjugate to $a\cdot S$, where the formal series $a$ has constant part 1.

\noindent
There exists $\zeta_0>0$ depending only on $\Lambda, \omega, n,d$ such that if $\lvert \lvert R\rvert\rvert_{r_0,\delta_0}\leq \zeta_0$, then $S+R$ is holomorphically conjugate to a resonant vector field $b\cdot S$, where $b$ is a holomorphic function.
\end{thm}


\subsubsection{Sequences of parameters}

Let $(m_k)_{k\geq 0}$ and $(r_k)_{k\geq 0}$ be the sequences defined in the assumption \ref{gamma+omega} (by assumption \ref{gamma+omega}.\ref{poslimitrk}, the sequence $(r_k)$ has a positive limit as $k\rightarrow \infty$). For all $k\geq 0$, let

\begin{equation}
\delta_{k+1}=\delta_k g(m_k)^{-(17+10n)/m_k}
\end{equation}

\noindent (recall that $m_k\geq 1$ and that $g(m_k)$ is positive). By the assumption \ref{gamma+omega}. \ref{brunosum}, the sequence $\delta_k$ has a positive limit $\delta_\infty\leq 1$ as $k\rightarrow +\infty$. 

\bigskip
\noindent \textbf{Notation:} For all $k\in\mathbb{N}$, we will hereafter denote by $\mathcal{O}_k$ the set $\mathcal{V}_{r_k}\times \mathcal{W}_{\delta_k}$.

\bigskip
Let $(\zeta_k)_{k\geq 0}$ be a real sequence such that

\begin{equation}\label{zeta0}\zeta_0<\frac{\delta_\infty}{8(n+d)(1+2C''_S (n+2)B)}
\end{equation}

\noindent and 

$$\zeta_{k+1}= \frac{2C''_S}{ g(m_k)}\zeta_k $$

\noindent where $C''_S$ is a fixed constant, depending only on $n,d$ and $S$, which will appear in the computations below, and $B$ is the Brjuno sum defined in the assumption \ref{gamma+omega}.\ref{brunosum}.
As mentioned in section \ref{setup}, one can assume that for all $k\geq 0$,

\begin{equation}\label{zetakdecr} \frac{2C''_S}{ g(m_k)}<1\end{equation}

\noindent so that $(\zeta_k)_{k\geq 0}$ is a strictly decreasing sequence.
Moreover, let 

$$\eta_k=\sum_{j=0}^k \zeta_j$$

\noindent Notice that for all $k\geq 0$, 

$$\eta_k=\zeta_0+\sum_{j=1}^k \zeta_0 \prod_{l=0}^{j-1} \frac{2C''_S}{g(m_l)}\leq \zeta_0+ \sum_{j=1}^k \zeta_0  \frac{2C''_S}{g(m_{j-1})}$$

\noindent Using the assumption \ref{gamma+omega}.\ref{convergence}, this implies that 

$$\eta_k\leq \zeta_0+\sum_{j=1}^k \zeta_0 \frac{2C''_S (n+2) \ln g(m_{j-1})}{m_{j-1}}\leq \zeta_0(1+2C''_S (n+2)B)$$

\noindent therefore, under the assumption \eqref{zeta0},
$\eta_k\leq \frac{1}{8}$.

%

\subsubsection{Iteration step}

\bigskip
\noindent
Theorem \ref{mainthmbis} will be proved by an iteration of the following statement, the proof of which is postponed to the last section:

\begin{prop}\label{iteration}
Let $N_k,R_k\in VF_{r_k,\delta_k}$ with $N_k$ resonant and $R_k$ of order $m_k$. Assume that 
\begin{enumerate}
\item $S+N_k+R_k$ is formally conjugate to the normal form $a\cdot S$ where $a$ is a formal series;
\item $\lvert\lvert N_k-S\rvert\rvert_{r_k,\delta_k} \leq \eta_k$,
\item \label{smallRk} $\lvert\lvert R_k\rvert\rvert_{r_k,\delta_k} \leq \zeta_k$,
\end{enumerate}
then there exists an analytic diffeomorphism $\Phi_k$ on 
$\mathcal{O}_{k+1}$
 with values in 
 $\mathcal{O}_k$, conjugating $S+N_k+R_k$ to $S+N_{k+1}+R_{k+1}$ where $N_{k+1}$ is resonant, $R_{k+1}$ has order $m_{k+1}$, such that the following properties hold:

\begin{enumerate}
\item \label{1} $\lvert \lvert N_{k+1}-S\rvert \rvert_{r_{k+1},\delta_{k+1}}\leq \eta_{k+1}$,
\item \label{2} $\lvert\lvert  R_{k+1}\rvert \rvert_{r_{k+1},\delta_{k+1}} \leq  \zeta_{k+1}$,
\item \label{3} $\lvert \lvert \Phi_k-Id\rvert\rvert_{r_{k+1},\delta_{k+1}}\leq \zeta_{k+1}$,
\item \label{DPhik} the quantity

$$\lvert \lvert D\Phi_k-I\rvert \rvert_{r_{k+1},\delta_{k+1}}\\
:= 
\max\{ \lvert \lvert \partial_{X_j}(\Phi_k-Id)\rvert \rvert_{r_{k+1},\delta_{k+1}}, \lvert \lvert \partial_{Y_{j'}}(\Phi_k-Id)\rvert \rvert_{r_{k+1},\delta_{k+1}}\}$$

\noindent is less than $ \zeta_{k+1}$.
\end{enumerate}

\end{prop}

\rem 
\noindent Under the assumptions of Theorem \ref{mainthmbis}, the initial system satisfies the assumptions of Proposition \ref{iteration} with $k=0$, $N_0=0$, $R_0=R$.

\subsubsection{Convergence of the algorithm}

\bigskip

\noindent We now proceed with the iteration of Proposition \ref{iteration}. 
For all $k\in\mathbb{N}$, let 

\begin{equation}\label{Psik}\Psi_k:=\Phi_0 \circ \dots \circ \Phi_k
\end{equation}

The map $\Psi_k$ is well defined on $\mathcal{O}_{k+1}$ with values in $\mathcal{O}_0$;
moreover, $\Psi_k$ conjugates $S+N_k+R_k$ to $S+R$. 

\begin{lem}
The map $\Psi_k$ is a holomorphic diffeomorphism from 
$\mathcal{O}_{k+1}$ to $\Psi_k(\mathcal{O}_{k+1})$.
\end{lem}

\dem The map $\Psi_k$ 
is analytic on $\mathcal{O}_{k+1}$ by composition, since every $\Phi_j$ is analytic on $\mathcal{O}_{j+1}$.
Moreover, for all $k$ and all $(X,Y)\in\mathcal{O}_{k+1}$,

\begin{equation}\label{Psik-Id}\begin{split}\lvert\lvert  D\Psi_k(X,Y)-I\rvert \rvert &
= \lvert \lvert( D\Psi_{k-1}-I)\circ \Phi_k \cdot D\Phi_k +D\Phi_k-I \rvert \rvert \\
&\leq \lvert \lvert( D\Psi_{k-1}-I)\circ \Phi_k \cdot D\Phi_k\rvert\rvert + \zeta_{k+1}
\end{split}\end{equation}

\noindent which implies, by the property \ref{DPhik} of Proposition \ref{iteration}, that on $\mathcal{O}_{k+1}$,

\begin{equation}\label{dpsikzetak}\lvert\lvert  D\Psi_k(X,Y)-I\rvert \rvert  \leq  \lvert \lvert( D\Psi_{k-1}-I)\circ \Phi_k \rvert\rvert (1+\zeta_{k+1}) + \zeta_{k+1}
\end{equation}

\noindent
Let $(u_n)_{n\in\mathbb{N}}$ be defined by $u_0=\zeta_1$ and for all $n\geq 1$, 

$$u_n=(1+\zeta_n)u_{n-1}+\zeta_n$$

\noindent One computes easily that for all $n\geq 2$,

$$u_n=\zeta_1\prod_{l=1}^n(1+\zeta_l)+\zeta_1\prod_{l=2}^n(1+\zeta_l)+\dots + \zeta_{n-2}(1+\zeta_{n-1})(1+\zeta_{n})+ \zeta_{n-1}(1+\zeta_n)+\zeta_n$$

\noindent therefore $u_n\leq \eta_n\prod_{l=1}^n(1+\zeta_l)\leq \frac{1}{8}e^{\eta_n}\leq \frac{1}{8}e^{\frac{1}{8}}$. Now, assuming
for all $(X,Y)\in \mathcal{O}_k$ that 

$$\lvert \lvert( D\Psi_{k-1}-I)(X,Y) \rvert\rvert \leq u_k$$

\noindent (which is true for $k=1$ by the property \ref{DPhik} of Proposition \ref{iteration} with $k=0$),
then on $\mathcal{O}_{k+1}$, the estimate \eqref{dpsikzetak} implies that

\begin{equation}\lvert\lvert  D\Psi_k(X,Y)-I\rvert \rvert  \leq u_k (1+\zeta_{k+1}) + \zeta_{k+1}=u_{k+1}
\end{equation}

\noindent
and
one deduces by recurrence that for all $k$, 

\begin{equation}\label{estimPsik}\lvert\lvert  D\Psi_k(X,Y)-I\rvert \rvert\leq \frac{1}{8}e^{\frac{1}{8}}<1
\end{equation}

\noindent
The map $\Psi_k$ is therefore an analytic diffeomorphism from $\mathcal{O}_{k+1}$ to its image. $\Box$

\begin{definition}
Let $\mathcal{O}_\infty=\mathcal{V}_{r_\infty}\times \mathcal{W}_{\delta_\infty}$, where $r_\infty$ is the limit of the sequence $(r_k)$ and $\delta_\infty$ is the limit of the sequence $(\delta_k)$, or equivalently, $\mathcal{O}_\infty=\cap_k \mathcal{O}_k$.
\end{definition}

%
%
%
%
%
%
%
%

\begin{prop}
The sequence $(\Psi_k)_{k\geq 1}$ defined in \eqref{Psik} has a subsequence which converges to a holomorphic diffeomorphism on $\mathcal{O}_\infty:=\cap_{k\in\mathbb{N}}\mathcal{O}_k$.
\end{prop}

\dem 
A simple recurrence shows that $\lvert \lvert Id-\Psi_k\rvert \rvert _{r_\infty,\delta_\infty}\leq \eta_{k+1} $. Indeed, 

$$Id - \Psi_k = (Id-\Psi_{k-1} )\circ \Phi_k +(Id-\Phi_k)$$

\noindent therefore

$$\lvert \lvert Id-\Psi_k\rvert \rvert _{r_{k+1},\delta_{k+1}} \leq \lvert  \lvert Id- \Psi_{k-1}\rvert \rvert _{r_k,\delta_k} +\lvert \lvert Id -\Phi_k\rvert \rvert _{r_{k+1},\delta_{k+1}}$$

\noindent and one sees that if $\lvert\lvert  Id- \Psi_{k-1}\rvert \rvert _{r_k,\delta_k}\leq \eta_k$, then $\lvert\lvert  Id-\Psi_k\rvert \rvert _{r_{k+1},\delta_{k+1}} \leq \eta_{k+1}$.
Thus $(\Psi_k)$ is uniformly bounded on $\mathcal{O}_\infty$ in the analytic weighted norm, therefore there exists a subsequence $(\Psi_{k_l})$ which is a Cauchy sequence for this norm. In particular, $(\Psi_{k_l})$ converges uniformly on every compact subset of $\mathcal{O}_\infty$ to a map $\Psi_\infty$ which is holomorphic on $\mathcal{O}_\infty$.

\bigskip
\noindent
Moreover, the estimate \eqref{estimPsik} is uniform in $k$. 
%
%
Therefore, $\Psi_\infty$ is still injective on $\mathcal{O}_\infty$, thus it is a holomorphic diffeomorphism on $\mathcal{O}_\infty$.
$\Box$

\bigskip
\noindent 
The sequence $(N_k)=(a_kS)$ is convergent in the topology of the analytic functions on $\mathcal{O}_\infty$, since $\lvert \lvert N_{k+1}-N_k\rvert \rvert_{r_\infty,\delta_\infty} =\lvert \lvert \bar{R}_{k,res}\rvert \rvert_{r_\infty,\delta_\infty}\leq \zeta_k$ and since $\zeta_k$ is summable (there exists a constant $c<1$ such that $\zeta_k\leq c^k \zeta_0$). Let $N_\infty$ be the limit of $(N_k)$ in this topology; then $N_\infty=a_\infty S$ for some resonant function $a_\infty$ which is holomorphic on $\mathcal{O}_\infty$. The sequence $R_k$ tends to 0 in the analytic topology on $\mathcal{O}_\infty$. 
Therefore $\Psi_\infty$ conjugates $S+N_\infty$ to $S+R$. 

\bigskip
\noindent
This concludes the proof of Theorem \ref{mainthmbis} based on Proposition \ref{iteration}, which is proved in the next section.

%
%
%
%
%
%
%
%

\section{Proof of Proposition \ref{iteration}}

We now make all the assumptions of Proposition \ref{iteration}.

\bigskip
\noindent 
If $S+N_k+R_k$ is formally conjugate to the normal form $a\cdot S$, then by Proposition \ref{nfk=aS}, there exists a formal series $a_k$ with constant part 1 such that $N_k=a_k\cdot S$ and since $N_k$ is analytic on $\mathcal{V}_{r_k}\times \mathcal{W}_{\delta_k}$, so is $a_k$. Moreover, if $\lvert\lvert  N_k-S\rvert \rvert_{r_k,\delta_k} \leq \eta_k$, then by our definition of the norm of a vector field given in section \ref{Notations}, for all $1\leq j\leq d$,

\begin{equation}
\lvert (a_k-1) \omega_j\rvert_{r_k,\delta_k} \leq \eta_k
\end{equation}

\noindent therefore $\lvert a_k-1\rvert_{r_k,\delta_k} \leq \frac{\eta_k}{\max_j \lvert \omega_j\rvert }\leq \frac{\delta_\infty}{8(n+d) C_\omega}\leq \frac{1}{2}$ (since it was assumed that $\delta_0\leq C_\omega$). 

%
%
%
%
\subsection{Homological equation}\label{homoleq}

As a first step in the construction of the diffeomorphism $\Phi_k$, one has to solve the homological equation 

$$\mathcal{L}G_k:=[G_k,N_k]=\bar{R}_{k}$$

\noindent 
where $\bar{R}_{k}$ contains all non resonant monomials of $R_k$ with degree between $m_k$ and $m_{k+1}-1$. 
For every vector field $F$, one has

$$\mathcal{L}F=[F,a_k S]=a_k[F,S]+F(a_k) S$$

\noindent
(where $F(a_k)$ stands for the Lie derivative of $a_k$ under $F$). Let $\mathcal{D}: F\mapsto a_k[F,S]$ and $\mathcal{N}:F\mapsto F(a_k)S$ 
which are defined and linear on the space of formal vector fields.
%

\noindent The operator $\mathcal{N}$ satisfies:

$$\mathcal{N}^2 F=\mathcal{N}(F(a_k)S)=F(a_k)S(a_k)=0$$

\noindent 
(since $a_k$ is resonant). Moreover, 

$$\mathcal{DN} F= \mathcal{D}(F(a_k)S)= a_k[F(a_k)S,S]=-a_k S(F(a_k))S$$

\noindent 
and 

$$\mathcal{ND} F= \mathcal{N} (a_k [F,S]) =-a_k S(F(a_k))S$$

\noindent 
so that $\mathcal{N} $ and $\mathcal{D}$ commute. 
Thus, formally

$$(\mathcal{D}+\mathcal{N})^{-1}=(I+\mathcal{D}^{-1}\mathcal{N})^{-1}\mathcal{D}^{-1}=\mathcal{D}^{-1}(I-\mathcal{ND}^{-1})$$

\noindent 
where $\mathcal{D}^{-1}$ and $(\mathcal{D}^{-1})^2$ are defined, i.e on non resonant formal vector fields.
Therefore, the solution will be 

$$G_k=\mathcal{D}^{-1}(I-\mathcal{ND}^{-1})\bar{R}_{k}$$

%



\noindent 
since $\bar{R}_k$ is a non resonant vector field.
Now for all $(P,Q)\in\mathbb{Z}^d\times \mathbb{N}^n$,

\begin{equation}\begin{split}ad_S&\left(\sum_{j=1}^d e^{i\langle P,X\rangle}Y^Q\frac{\partial}{\partial X_j} + \sum_{j'=1}^n e^{i\langle P,X\rangle}Y^{Q+E_{j'}}\frac{\partial}{\partial Y_{j'}}\right)\\
&=(i\langle P,\omega\rangle +\langle Q,\Lambda\rangle )\left(\sum_{j=1}^d e^{i\langle P,X\rangle}Y^Q\frac{\partial}{\partial X_j} + \sum_{j'=1}^n e^{i\langle P,X\rangle}Y^{Q+E_{j'}}\frac{\partial}{\partial Y_{j'}}\right)
\end{split}\end{equation}

\noindent
therefore the operators $ad_S$ and $a_k\cdot$ (the multiplication by $a_k$) preserve the equivalence classes $\mathcal{C}_c$ 
of indices $(P,Q)\in \mathbb{Z}^d\times \mathbb{N}^n$ such that $ i\langle P,\omega\rangle +\langle Q,\Lambda\rangle=c$;
moreover, $\mathcal{N}$ also preserves these equivalence classes since 

\begin{equation}\begin{split}&\sum_{j=1}^d e^{i\langle P,X\rangle}Y^Q\frac{\partial}{\partial X_j}+ \sum_{j'=1}^n e^{i\langle P,X\rangle}Y^{Q+E_{j'}}\frac{\partial}{\partial Y_{j'}}\\
&=\left( \sum_{j=1}^d e^{i\langle P,X\rangle}Y^Q\frac{d}{d X_j} a_k(X,Y)+ \sum_{j'=1}^n e^{i\langle P,X\rangle}Y^{Q+E_{j'}}\frac{d}{d Y_{j'}} a_k(X,Y)
\right) S(Y)\\
\end{split}
\end{equation}

Therefore,

\begin{equation}\begin{split}&(\mathcal{D}+\mathcal{N})F=R \Leftrightarrow F=\sum_c F_c \ \mathrm{and}\ \forall c\in \mathbb{C}, \\
& (\mathcal{D}+\mathcal{N}) F_c = 
\sum_{(P,Q)\in\mathcal{C}_c} \left(\sum_{j=1}^d R_{j,P,Q} \frac{\partial}{\partial X_j}+\sum_{j'=1}^nR_{d+j',P,Q+E_{j'}}Y_{j'}\frac{\partial}{\partial Y_{j'}}\right) e^{i\langle P,X\rangle } Y^Q
\end{split}\end{equation}

\noindent 
This makes it possible to separate the homological equation into two parts: the low frequency part, and the high frequency part. This means distinguishing two cases for $c$:
\begin{enumerate}
\item 
either $\mathcal{C}_c$ contains a couple $(P,Q)$ with $\lvert P\rvert \leq m_k$, 
\item or all elements of $\mathcal{C}_c$ satisfy $\lvert P\rvert >m_k$ 
\end{enumerate}
(the equivalence class $\mathcal{C}_0$ does not enter into the computation, since $\bar{R}_k$ is non resonant). We will denote by $\mathcal{I}_0$ the set of all non empty equivalence classes which are in the first case, and by $\mathcal{I}_\infty$ the set of all equivalence classes which are in the second case.
Thus there is a decomposition 

$$\bar{R}_k=\bar{R}_k^0 +\bar{R}_k^\infty$$

\noindent where $\bar{R}_k^0$ only has monomials $\bar{R}_{k,P,Q}$ such that $(P,Q)$ belongs to an element of $\mathcal{I}_0$, and $\bar{R}_k^\infty$ only has monomials $\bar{R}_{k,P,Q}$ such that $(P,Q)$ belongs to an element of $\mathcal{I}_\infty$, both being non resonant.
Since $\bar{R}_k$ is truncated with respect to the degree in $Y$, the low frequency part $\bar{R}_k^0$ has monomials whose indices belong to a finite number of equivalence classes.

\bigskip

\noindent
One will solve separately $\mathcal{L}G_k^0=\bar{R}_k^0$ and $\mathcal{L} G_k^\infty=\bar{R}_k^\infty$, then let $G_k=G_k^0+G_k^\infty$. As shown above, the solutions are given by

$$G_k^0=\mathcal{D}^{-1}(I-\mathcal{ND}^{-1})\bar{R}_{k}^0
$$

\noindent and

$$G_k^\infty=\mathcal{D}^{-1}(I-\mathcal{ND}^{-1})\bar{R}_{k}^\infty
$$

\noindent (since $\bar{R}_{k}^0$ and $\bar{R}_{k}^\infty$ are non resonant, the solutions are defined and they will be analytic under some assumptions on $\bar{R}_k^0,\bar{R}_k^\infty$). It remains to estimate $G_k^0$ and $G_k^\infty$ in the analytic norm on $\mathcal{O}_{k+1}$.

\paragraph{Estimate of the operator $\mathcal{N}$}

\noindent
It is convenient first to give an estimate of $\mathcal{N}$ acting on $VF_{r'',\delta''}$ into $VF_{r_3,\delta_3}$, where the positive parameters $r''>r_3,\delta''>\delta_3$ will be given later. Applying $\mathcal{N}$ does not add small divisors.
One has

\begin{equation}
\begin{split}F(a_k)(X,Y)
&=\sum_{j=1}^d F_j(X,Y) \frac{d}{dX_j} a_k(X,Y) +\sum_{j'=1}^n F'_{j'} (X,Y) \frac{d}{dY_{j'}} a_k(X,Y)\\
&=\sum_{P,Q} (a_k)_{P,Q} \sum_{j=1}^d F_j(X,Y) iP_j e^{i\langle P,X\rangle }Y^Q \\
&+\sum_{P,Q/ Q_{j'}\geq 1} (a_k)_{P,Q} \sum_{j'=1}^n F'_{j'} (X,Y) Q_{j'} Y^{Q-E_{j'}}e^{i\langle P,X\rangle } \\
\end{split}
\end{equation}

%

\noindent which implies that 

\begin{equation}\label{estimateN}\begin{split}&\lvert F(a_k)\rvert_{r_3,\delta_3}\\
&\leq \sum_{P,Q} e^{r_3\lvert P\rvert }Ê\delta_3^{\lvert Q\rvert } \sum_{P',Q'} \lvert (a_k)_{P',Q'}\rvert[\sum_j \lvert P'_j\rvert \lvert F_{j,P-P',Q-Q'}\rvert +\sum_{j'} \lvert Q'_{j'}\rvert \lvert F_{j',P-P',
Q+E_{j'}-Q'}\rvert ]\\
&\leq \lvert a_k\rvert_{r'',\delta''} \lvert\lvert  F\rvert\rvert_{r'',\delta''}\sum_{P,Q} e^{r_3\lvert P\rvert }Ê\delta_3^{\lvert Q\rvert } \sum_{P',Q'} e^{-\lvert P'\rvert r''} \delta''^{-\lvert Q'\rvert}e^{-\lvert P-P'\rvert r''} \delta''^{-\lvert Q-Q'\rvert } [\sum_j \lvert P'_j\rvert \\
&+\sum_{j'} \lvert Q'_{j'}\rvert \delta'']\\
\end{split}
\end{equation}

\noindent Let 
$r_4=\frac{r_3+r''}{2}$ and $\delta_4=\sqrt{\delta'' \delta_3}$, in order to have 
$r''-r_4=r_4-r_3$ and $\delta''/\delta_4=\delta_4/\delta_3$.
Then

\begin{equation}\begin{split}&\lvert F(a_k)\rvert_{r_3,\delta_3}\\
&\leq 
\lvert a_k\rvert_{r'',\delta''} \lvert\lvert  F\rvert\rvert_{r'',\delta''}\sum_{P,Q} e^{r_3\lvert P\rvert }Ê\delta_3^{\lvert Q\rvert } \sum_{P'} e^{-\lvert P'\rvert r_4} e^{-\lvert P-P'\rvert r''}  e^{-\lvert P'\rvert (r''-r_4)}(\lvert P'\rvert +n)\\
&\qquad \qquad \qquad  \qquad\cdot \sum_{Q'}\delta_4^{-\lvert Q'\rvert}( \delta''/\delta_4)^{-\lvert Q'\rvert}\delta''^{-\lvert Q-Q'\rvert }
(d+\lvert Q'\rvert \delta'')\\
&\leq 
\lvert a_k\rvert_{r'',\delta''} \lvert\lvert  F\rvert\rvert_{r'',\delta''}\sum_{P,Q} e^{r_3\lvert P\rvert }Ê\delta_3^{\lvert Q\rvert } \sum_{P'} e^{-\lvert P\rvert r_4}   e^{-\lvert P'\rvert (r''-r_4)}(\lvert P'\rvert +n)\\
&\qquad\qquad\qquad\qquad\cdot \sum_{Q'}\delta_4^{-\lvert Q\rvert }( \delta''/\delta_4)^{-\lvert Q'\rvert}
(d+\lvert Q'\rvert \delta'')\\
\end{split}
\end{equation}

\noindent Now 

\begin{equation}\label{sumfreq}\begin{split}
 \sum_{P'}  e^{-\lvert P'\rvert (r''-r_4)}(\lvert P'\rvert +n)&=\sum_{K\in\mathbb{N}} \sum_{\lvert P'\rvert=K} (K+n)e^{-K(r''-r_4)}\\
 &\leq \sum_{K\in\mathbb{N}} C(d) K^{d-1} (K+n) e^{-K(r''-r_4)}\\
 &\leq C(n,d) \int_0^\infty t^d e^{-t(r''-r_4)}dt\\
 &\leq \frac{C(n,d)}{(r''-r_4)^{d+1}}\\
\end{split}\end{equation}

\noindent and 

\begin{equation}\label{sumpuiss}\begin{split}
\sum_{Q'}( \delta''/\delta_4)^{-\lvert Q'\rvert}
(d+\lvert Q'\rvert \delta'') 
&= \sum_{L\in\mathbb{N}} \sum_{\lvert Q'\rvert =L} ( \delta''/\delta_4)^{-L}
(d+L\delta'') \\
&\leq C(n,d)\sum_{L\in\mathbb{N}} L^n ( \delta''/\delta_4)^{-L}\\
&\leq C(n,d)\int_0^\infty t^n ( \delta''/\delta_4)^{-t}dt\\
&\leq C(n,d)(\ln \delta ''/\delta_4)^{-n-1}\\
\end{split}\end{equation}

\noindent (here and below, $C(d) $ and $C(n,d)$ stand for generic constants which depend only on $n,d$) therefore

\begin{equation}\begin{split}&\lvert F(a_k)\rvert_{r_3,\delta_3}\\
&\leq C(n,d) \lvert a_k\rvert_{r'',\delta''} \lvert \lvert F\rvert\rvert_{r'',\delta''}\sum_{P,Q} e^{r_3\lvert P\rvert }Ê\delta_3^{\lvert Q\rvert } e^{-\lvert P\rvert r_4} (r''-r_4)^{-(d+1)}\delta_4^{-\lvert Q\rvert } (\ln \delta ''/\delta_4)^{-n-1}\\
\end{split}
\end{equation}

\noindent
Thus, 

$$\lvert F(a_k)\rvert_{r_3,\delta_3}\leq \lvert\lvert  F\rvert\rvert_{r'',\delta''} \lvert a_k\rvert_{r'',\delta''} \frac{C(n,d) }{ (r''-r_3)^{2d+4}(\ln \delta_3-\ln \delta'')^{2n+4}}$$

\noindent
Therefore, 

\begin{equation}\label{opN}\lvert \lvert \lvert \mathcal{N}\rvert\rvert \rvert_{C^\omega_{r'',\delta''}\rightarrow C^\omega_{r_3,\delta_3}} \leq \lvert a_k\rvert_{r'',\delta''} \frac{C(n,d) }{ (r''-r_3)^{2d+4}(\ln \delta_3-\ln \delta'')^{2n+4}}C_\omega
\end{equation}

%


\paragraph{Low frequency part}

\noindent
First consider the set $\mathcal{I}_0$ and let $\mathcal{C}_c\in\mathcal{I}_0$ with $c\neq 0$. 
Let 

$$R_c(X,Y)=\sum_{(P,Q)\in \mathcal{C}_c} \left(\sum_{j=1}^d (\bar{R}_{k})_{j,P,Q} \frac{\partial}{\partial X_j} +\sum_{j'=1}^n (\bar{R}_k)_{d+j',P,Q+E_{j'}} Y_{j'}\frac{\partial}{\partial Y_{j'}}
\right)e^{i\langle P,X\rangle } Y^Q$$

\noindent The equation 
$\mathcal{D} F_c =R_c$ can be solved as follows:

\begin{equation}\label{equivclass}
F_c=\mathcal{D}^{-1}{R_c}=\sum_{(P,Q)\in\mathcal{C}_c} \frac{(R_ca_k^{-1})_{P,Q} }{c } e^{i\langle P,X\rangle } Y^Q
\end{equation}

\noindent
Now, since $a_k$ is analytically close to 1 (i.e. $\lvert a_k-1\rvert_{r_k,\delta_k} \leq \frac{1}{2}$), using the sub-multiplicativity of the weighted norm, $R_ca_k^{-1}$ is in $VF_{r_k,\delta_k}$, with norm less than $2\lvert\lvert R_c\rvert\rvert_{r_k,\delta_k}$ and $R_ca_k^{-1}$ also has quasi-order $m_k$ in $Y$. Thus, using the estimate \eqref{coeffPQ}, for all $1\leq j\leq d+n$, $(R_ca_k^{-1})_{j,P,Q}$ has modulus less than $2\lvert \lvert R_c\rvert\rvert_{r_k,\delta_k}e^{-r_k\lvert P\rvert} \delta_k^{-\lvert Q\rvert}$.  Therefore, letting $\delta=\delta_k,r=r_k$, for all $r'\in (0,r)$ and all $\delta'\in (0,\delta)$,

\begin{equation}\label{D-1}\begin{split}
\lvert\lvert  \mathcal{D}^{-1}{R_c}\rvert\rvert_{r',\delta'}&\leq\frac{2}{c}\lvert \lvert R_c\rvert \rvert_{r,\delta}\sum_{Q,\lvert Q\rvert= m_k}^{m_{k+1}-1}(\frac{\delta'}{\delta})^{\lvert Q\rvert} \sum_{P/(P,Q )\in\mathcal{C}_c}   e^{-(r-r')\lvert P\rvert} \\
&\leq \frac{2\lvert \lvert R_c\rvert\rvert_{r,\delta}}{c(r-r')^{d+1}}\sum_{Q,\lvert Q\rvert= m_k}^{m_{k+1}-1}(\frac{\delta'}{\delta})^{\lvert Q\rvert}\\
&\leq \frac{2\lvert \lvert R_c\rvert\rvert_{r,\delta}}{c(r-r')^{d+1}}(\frac{\delta '}{\delta})^{m_k}\frac{1}{(\ln \delta' -\ln \delta)^{n+1}}
\end{split}
\end{equation}

\noindent
This implies the following estimate:

\begin{equation}
\label{Dinverse}\lvert\lvert  \mathcal{D}^{-1}\bar{R}_k^0\rvert \rvert_{r',\delta'}\leq \frac{2\sup_{c\in \mathcal{I}_0\setminus\{0\}}\{c^{-1}\}}{(r-r')^{d+1}}(\ln \frac{\delta'}{\delta})^{-(n+1)}\lvert \lvert \bar{R}_k^0\rvert\rvert_{r,\delta}(\frac{\delta'}{\delta})^{m_k}
\end{equation}

\noindent (since $\lvert \lvert \bar{R}_k^0\rvert\rvert_{r,\delta}=\sum_{\mathcal{C}_c\in\mathcal{I}_0}\lvert \lvert R_c\rvert\rvert_{r,\delta}$).
Applying $\mathcal{D}^{-1}$ once more on \eqref{equivclass}, one has 

$$(\mathcal{D}^{-1})^2 R_c =\sum_{(P,Q)\in\mathcal{C}_c} \frac{((\mathcal{D}^{-1}R_c)a_k^{-1})_{P,Q} }{c } e^{i\langle P,X\rangle } Y^Q$$

\noindent
and since $(\mathcal{D}^{-1})R_c$ is in 
$C^\omega_{r',\delta'}$ and 
$\lvert a_k-1\rvert_{r',\delta'}\leq \frac{1}{2}$,
then for all $1\leq j\leq d+n$,

$$ \lvert ((\mathcal{D}^{-1}R_c)a_k^{-1})_{j,P,Q}  \rvert \leq \lvert\lvert  (\mathcal{D}^{-1}R_c)a_k^{-1}\rvert\rvert _{r',\delta'} e^{-\lvert P\rvert r'} \delta'^{-\lvert Q\rvert}$$

\noindent
thus for all $r''<r'$ and $\delta''<\delta'$,

\begin{equation}\begin{split}\lvert \lvert (\mathcal{D}^{-1})^2{R_c}\rvert\rvert_{r'',\delta''}&\leq \frac{4\sup_{c\in \mathcal{I}_0\setminus\{0\}}\{c^{-1}\}^2}{(r-r')^{d+1}(r'-r'')^{d+1}}(\frac{\delta''}{\delta})^{m_k}\\
&\cdot \frac{1}{(\ln \delta'-\ln \delta)^{n+1}(\ln \delta''-\ln \delta')^{n+1}}\lvert \lvert R_c\rvert \rvert_{r,\delta}
\end{split}\end{equation}

\noindent therefore, by the assumption \ref{gamma+omega}.\ref{lowfq},

\begin{equation}\begin{split}\lvert \lvert (\mathcal{D}^{-1})^2{R_c}\rvert \rvert_{r'',\delta''}&\leq \frac{4g(m_k)^2}{(r-r')^{d+1}(r'-r'')^{d+1}}(\frac{\delta''}{\delta})^{m_k}\\
&\cdot \frac{1}{(\ln \delta'-\ln \delta)^{n+1}(\ln \delta''-\ln \delta')^{n+1}}\lvert\lvert  R_c\rvert\rvert_{r,\delta}
\end{split}\end{equation}

Using now the estimate \eqref{opN}, for the low frequencies of $\bar{R}_{k}$, one has, for some constant $C'$ which only depends on $n,d$, 

\begin{equation}\begin{split}&\lvert\lvert  \mathcal{D}^{-1}(I-\mathcal{N}\mathcal{D}^{-1}) R_c \rvert\rvert_{r_3,\delta_3}\\
&\leq 
\frac{C'g(m_k)^2}{[(r-r')(r'-r'')]^{d+1}}\frac{(\frac{\delta''}{\delta})^{m_k}\lvert a_k\rvert_{r'',\delta''}}{[(\ln (\delta'/\delta))(\ln (\delta''/ \delta'))]^{n+1}} \frac{ C_\omega\lvert\lvert  R_c\rvert\rvert_{r,\delta}}{ (r''-r_3)^{2d+4}(\ln( \delta_3/\delta''))^{2n+4}}
\end{split}\end{equation}

\noindent
therefore

\begin{equation}\label{Gk0}\begin{split}
&\lvert\lvert  G_k^0\rvert \rvert_{r_3,\delta_3}\\
&\leq \frac{C'g(m_k)^2}{[(r-r')(r'-r'')]^{d+1}}\frac{(\frac{\delta''}{\delta})^{m_k}\lvert a_k\rvert_{r'',\delta''}}{[(\ln (\delta'/\delta))(\ln (\delta''/ \delta'))]^{n+1}} \frac{C_\omega\lvert \lvert \bar{R}_k^0\rvert\rvert_{r,\delta}}{ (r''-r_3)^{2d+4}(\ln( \delta_3/\delta''))^{2n+4}}
\end{split}\end{equation}

\paragraph{High frequency part}

If $\mathcal{C}_c$ is in the second case (which implies that $c\neq 0$), by the Assumption \ref{gamma+omega}, one has the estimate

\begin{equation}
\begin{split}\lvert\lvert  \mathcal{D}^{-1}R_c\rvert \rvert_{r',\delta'}&\leq 
\frac{2}{c}\lvert\lvert  R_c\rvert \rvert_{r,\delta}\sum_{Q,\lvert Q\rvert\geq m_k}(\frac{\delta'}{\delta})^{\lvert Q\rvert} \sum_{P/(P,Q )\in\mathcal{C}_c}   e^{-(r-r')\lvert P\rvert} \\
&\leq 2\lvert \lvert R_c\rvert\rvert_{r,\delta}\sum_{Q,\lvert Q\rvert\geq m_k}(\frac{\delta'}{\delta})^{\lvert Q\rvert} \sum_{P/(P,Q )\in\mathcal{C}_c}   e^{-(r-r'-\epsilon_k)\lvert P\rvert} \\
\end{split}
\end{equation}

\noindent (the computation is similar to \eqref{D-1})
thus

\begin{equation}\begin{split}
\lvert \lvert \mathcal{D}^{-1}\bar{R}_k^\infty\rvert \rvert_{r',\delta'}&\leq 2\sum_{\mathcal{C}_c\in \mathcal{I}_\infty} \lvert \lvert R_c\rvert\rvert_{r,\delta}\sum_{Q,\lvert Q\rvert\geq m_k}(\frac{\delta'}{\delta})^{\lvert Q\rvert} \sum_{P/(P,Q )\in\mathcal{C}_c}   e^{-(r-r'-\epsilon_k)\lvert P\rvert} \\
&\leq 2 \lvert\lvert  \bar{R}_k^\infty \rvert\rvert_{r,\delta}\sum_{Q,\lvert Q\rvert\geq m_k}(\frac{\delta'}{\delta})^{\lvert Q\rvert} \sum_{P, \lvert P\rvert > m_k}   e^{-(r-r'-\epsilon_k)\lvert P\rvert} \\
&\leq  \frac{2 \lvert\lvert  \bar{R}_k^\infty\rvert \rvert_{r,\delta}}{(r-r'-\epsilon_k)^{d+1}}(\frac{\delta '}{\delta})^{m_k}\frac{e^{-(r-r'-\epsilon_k)m_k}}{(\ln \delta' -\ln \delta)^{n+1}}
\end{split}\end{equation}

\noindent (as in \eqref{sumfreq} and \eqref{sumpuiss}, the sums on $P$ and $Q$ were estimated by means of an integral).
Applying $\mathcal{D}^{-1}$ once more, one gets

\begin{equation}
\begin{split}
&\lvert\lvert ( \mathcal{D}^{-1})^2\bar{R}_k^\infty\rvert \rvert_{r'',\delta''}\\
& \leq \frac{4 \lvert\lvert  \bar{R}_k^\infty \rvert \rvert_{r,\delta}}{[(r-r'-\epsilon_k)(r'-r''-\epsilon_k)]^{d+1}}(\frac{\delta ''}{\delta})^{m_k}\frac{e^{-(r-r'-\epsilon_k)m_k}e^{-(r'-r''-\epsilon_k)m_k}}{[(\ln (\delta' / \delta)\ln (\delta''/\delta')]^{n+1}}\\
&\leq  \frac{4 \lvert \lvert \bar{R}_k^\infty \rvert\rvert_{r,\delta}}{[(r-r'-\epsilon_k)(r'-r''-\epsilon_k)]^{d+1}}(\frac{\delta ''}{\delta})^{m_k}\frac{e^{-(r-r''-2\epsilon_k)m_k}}{[(\ln (\delta' / \delta)\ln (\delta''/\delta')]^{n+1}}\\
\end{split}
\end{equation}

\noindent
Finally, using the estimate \eqref{opN},

\begin{equation}\label{Rkinfty}\begin{split}
&\lvert\lvert \mathcal{N} (\mathcal{D}^{-1})^2\bar{R}_k^\infty\rvert \rvert_{r_3,\delta_3}\\
&\leq  \frac{4 \lvert \lvert \bar{R}_k^\infty \rvert \rvert_{r,\delta}(\delta ''/\delta)^{m_k}}{[(r-r'-\epsilon_k)(r'-r''-\epsilon_k)]^{d+1}}\frac{e^{-(r-r''-2\epsilon_k)m_k}\lvert a_k\rvert_{r'',\delta''} }{[(\ln (\delta' / \delta)\ln (\delta''/\delta')]^{n+1}}\\
&\cdot \frac{C_{n,d}C_\omega}{ (r''-r_3)^{2d+4}(\ln \delta_3-\ln \delta'')^{2n+4}}
\end{split}\end{equation}

\noindent which implies that

\begin{equation}\label{Gkinfty}\begin{split}
\lvert\lvert  G_k^\infty\rvert \rvert_{r_3,\delta_3}
&\leq  \frac{C_{n,d,S} \lvert\lvert  \bar{R}_k^\infty\rvert \rvert_{r,\delta}(\delta ''/\delta)^{m_k}}{[(r-r'-\epsilon_k)(r'-r''-\epsilon_k)]^{d+1}}\frac{e^{-(r-r''-2\epsilon_k)m_k} }{[(\ln (\delta' / \delta)\ln (\delta''/\delta')]^{n+1}}\\
&\cdot [ (r''-r_3)^{2d+4}(\ln \delta_3-\ln \delta'')^{2n+4}]^{-1}
\end{split}\end{equation}

\noindent where $C_{n,d,S}$ depends on $n,d,S$ but not on the step of iteration.

\subsection{Choice of the parameters}

We make the following choice of parameters: 

$$\delta=\delta_k;\delta''=\delta g(m_k)^{\frac{-D''}{m_k}}; \ \delta'=\sqrt{\delta''\delta};\ \delta_3=\delta^{''2}/\delta'\ \mathrm{and} \ \delta_4=\delta_{k+1}$$

\noindent where $D''=10+6n$, and

\begin{equation}\label{width} r'=r - \frac{1}{4} (r_k-r_{k+1}); r''=r- \frac{1}{2} (r_k-r_{k+1});\ r_3=r- \frac{3}{4} (r_k-r_{k+1}); r_4=r_{k+1}
\end{equation}

\noindent 
This choice implies that $r-r'=r'-r''=r''-r_3$ and 

$$\lvert \ln (\delta'/\delta)\rvert = \lvert \ln (\delta''/\delta')\rvert =\lvert \ln (\delta_3/\delta'')\rvert =\frac{1}{2}\lvert \ln (\delta''/\delta)\rvert  = \frac{D''}{2m_k}\ln g(m_k)$$

\noindent as well as

$$\ln \delta_3-\ln \delta_4 = (n+2) \frac{\ln g(m_k)}{m_k}$$

\noindent 
Applied to estimate \eqref{Gk0} on the low frequency part $G_k^0$, this choice gives

\begin{equation}\label{Gklow}
\lvert\lvert  G_k^0\rvert \rvert_{r_3,\delta_3}
\leq \frac{C'\lvert a_k\rvert_{r'',\delta''}g(m_k)^2 (\frac{\delta''}{\delta})^{m_k}}{ (r-r')^{4d+6}[ \frac{D''}{2m_k}\ln g(m_k)  ]^{4n+6}} C_\omega\lvert\lvert  \bar{R}_k\rvert\rvert_{r,\delta}\end{equation}

\noindent 
Now by assumption \ref{gamma+omega}.\ref{rkassezgrand},

$$(r-r')^{-4d-6}= [\frac{1}{4}(r_k-r_{k+1})]^{-4d-6}\leq 4^{4d+6} (r_k-r_{k+1}-\epsilon_k)^{-4d-6}\leq g(m_k)$$

\noindent 
Therefore,

\begin{equation}
\lvert\lvert  G_k^0\rvert \rvert_{r_3,\delta_3}\leq \frac{C''_Sg(m_k)^3(\frac{\delta''}{\delta})^{m_k}}{[ \frac{D''}{2m_k}\ln g(m_k)  ]^{4n+6}} \lvert \lvert \bar{R}_k\rvert \rvert_{r,\delta}\end{equation}

\noindent where $C''_S$ only depends on $n,d,S$. 
By Assumption \ref{gamma+omega}. \ref{convergence},

\begin{equation}
\lvert \lvert G_k^0\rvert \rvert_{r_3,\delta_3}\leq C''_Sg(m_k)^{3+4n+6-D''}\lvert\lvert  \bar{R}_k\rvert \rvert_{r,\delta}
\leq \frac{ C'' _S}{ g(m_k)}\lvert \lvert \bar{R}_k\rvert\rvert_{r,\delta}\end{equation}

\noindent since $D''= 10+6n$. By the smallness assumption \ref{smallRk} on $\bar{R}_k$ and the assumption \eqref{zetakdecr} on $C''_S$ and $g$, and by the definition of the sequence $(\zeta_k)$,

\begin{equation}
\lvert \lvert G_k^0\rvert\rvert_{r_3,\delta_3}\leq 
 \frac{\zeta_{k+1}}{2}\end{equation}

\noindent
Concerning high frequencies, when applied to estimate \eqref{Gkinfty}, the choice of parameters implies

\begin{equation}\label{Gkhigh}
\lvert \lvert G_k^\infty\rvert \rvert_{r_3,\delta_3}
\leq  \frac{C_{n,d,S} \lvert \lvert \bar{R}_k^\infty\rvert  \rvert_{r,\delta}(\delta ''/\delta)^{m_k}}{(r-r'-\epsilon_k)^{2(d+1)}}\frac{e^{-(r-r''-2\epsilon_k)m_k} }{[ \frac{D''}{2m_k}\ln g(m_k)  ]^{4(n+1)+2n+4}}(\frac{1}{4}(r_k-r_{k+1}))^{-2d-4}
\end{equation}

\noindent
Again by Assumption \ref{gamma+omega}.\ref{convergence} and \ref{gamma+omega}.\ref{rkassezgrand}, this implies 

\begin{equation}
\lvert \lvert G_k^\infty\rvert \rvert_{r_3,\delta_3}
\leq   C_{n,d,S}\lvert\lvert  \bar{R}_k^\infty \rvert \rvert_{r,\delta}g(m_k)^{1+4(n+1)+2n+4}(\delta ''/\delta)^{m_k} e^{-(r-r''-2\epsilon)m_k}
\end{equation}

\noindent
Therefore, since $D''\geq 2+4(n+1)+2n+4$ and since $(\delta''/\delta)^{m_k}\leq g(m_k)^{-D''}$,

\begin{equation}
\lvert\lvert  G_k^\infty\rvert \rvert_{r_3,\delta_3}
\leq  \frac{ C_{n,d,S} }{ g(m_k)}\lvert \lvert \bar{R}_k^\infty \rvert \rvert_{r,\delta}\end{equation}

\noindent which implies, using the assumption \ref{smallRk} and the definition of the sequence $(\zeta_k)$, that 

$$\lvert \lvert G_k^\infty\rvert \rvert_{r_3,\delta_3}\leq \frac{\zeta_{k+1}}{2}$$

\noindent 
The following estimate of $G_k$ is finally obtained:

\begin{equation}\label{Gk}\lvert \lvert G_k\rvert \rvert_{r_3,\delta_3} \leq \zeta_{k+1}
\end{equation}

\bigskip
\noindent
Analytic bounds on the differential $DG_k$ of $G_k$ will also be necessary. For all $1\leq j\leq d$, 

\begin{equation}\begin{split}\lvert \lvert \partial_{X_j} G_{k}\rvert\rvert _{r_4,\delta_4}
&=\sum_{P,Q\geq m_k} \lvert \lvert iP_j G_{k,P,Q}\rvert \rvert e^{r_4\lvert P\rvert } \delta_4^Q\\
&\leq (\frac{\delta_4}{\delta_3})^{m_k} \lvert \lvert G_k\rvert \rvert _{r_3,\delta_3}(r_3-r_4)^{-(d+2)}(\ln \delta_3-\ln \delta_4)^{-(n+1)}
\end{split}\end{equation}

\noindent 
Similarly, for all $1\leq j'\leq n$, 

\begin{equation}\begin{split}\lvert \lvert \partial_{Y_{j'}} G_{k}\rvert\rvert _{r_4,\delta_4}
&=\sum_{P,Q\geq m_k} \lvert \lvert Q_{j'} G_{k,P,Q}\rvert \rvert e^{r_4\lvert P\rvert } \delta_4^{Q-1}\\
&\leq (\frac{\delta_4}{\delta_3})^{m_k}\lvert \lvert G_k\rvert \rvert _{r_3,\delta_3}(r_3-r_4)^{-(d+1)}(\ln \delta_3-\ln \delta_4)^{-(n+2)}
\end{split}\end{equation}

\noindent 
Therefore (using the estimate \eqref{Gk}, the assumption \ref{gamma+omega}.\ref{rkassezgrand} and the choice of $\delta_4$),

\begin{equation}\label{DGk}\begin{split}\lvert \lvert DG_k\rvert \rvert_{r_4,\delta_4}&:= \max_{j,j'}\{\lvert \lvert \partial_{X_j} G_{k}\rvert\rvert _{r_4,\delta_4}, \lvert \lvert \partial_{Y_{j'}} G_{k}\rvert\rvert _{r_4,\delta_4} \}\\
& \leq (\frac{\delta_4}{\delta_3})^{m_k}g(m_k)^{n+2} \zeta_{k+1}\leq \zeta_{k+1}
\end{split}\end{equation}




\subsection{Change of variables}



%
%

%
%

Let $G_k$ be the vector field on $\mathcal{O}_{{k+1}}$ with order $m_k+1$ defined in Section \ref{homoleq}, and let 
$\Phi_k=Id+G_k$ (this is a slight abuse of notation: we identify the vector $G_k(X,Y)\in\mathbb{C}^d\times \mathbb{C}^n$ with its projection on $\mathbb{C}^d/\mathbb{Z}^d \times \mathbb{C}^n$). 
Then 

$$\lvert\lvert  \Phi_k-Id\rvert \rvert _{r_{k+1},\delta_{k+1}}=\lvert \lvert G_k\rvert \rvert _{r_{k+1},\delta_{k+1}}\leq \zeta_{k+1}$$

\noindent therefore property \ref{3} holds. 
One sees that $\Phi_k$ is a diffeomorphism from 
$\mathcal{O}_{k+1}$ to $\mathcal{O}_k$: indeed, 
if $(X,Y)\in \mathcal{O}_{k+1}$, then for all $1\leq j\leq d$,

\begin{equation}
\lvert \Phi_{k,j}(X,Y)\rvert \leq \lvert X_j\rvert + \zeta_{k+1} \leq r_{k+1}+\zeta_{k+1}  
\end{equation}

\noindent and the assumption \ref{gamma+omega}.\ref{rkassezgrand} implies that $\lvert \Phi_{k,j}(X,Y)\rvert \leq r_k$, as long as $g(m_0)$ is large enough as a function of $d,S$, which is not a restrictive assumption. For all $1\leq j'\leq n$, 

\begin{equation}
\lvert \Phi_{k,j'}(X,Y)\rvert \leq \lvert Y_{j'}\rvert + \zeta_{k+1} \leq \delta_{k+1}+\zeta_{k+1}  
\end{equation}

\noindent therefore $\Phi_k$ has values in $\mathcal{O}_k$ if it can be shown that 

\begin{equation}\label{zetakdeltak}\frac{2C''_S\zeta_k}{g(m_k)}\leq \delta_k-\delta_{k+1}=\delta_{k+1}(g(m_k)^{\frac{17+10n}{m_k}}-1)
\end{equation}

\noindent Since by the assumption \eqref{zeta0}, $\zeta_0\leq \frac{1}{2C''_S}\delta_\infty$, then in particular $2C''_S\zeta_k\leq \delta_{k+1}$; moreover, 

$$g(m_k)^{\frac{17+10n}{m_k}}-1\geq \frac{17+10n}{m_k}\ln g(m_k)\geq  \frac{17+10n}{n+2}  g(m_k)^{-1}\geq g(m_k)^{-1}$$

\noindent which implies \eqref{zetakdeltak}.
Therefore $\Phi_k$ has values in $\mathcal{O}_k$.

\bigskip
\noindent
By the estimate \eqref{DGk}, for all $(X,Y)\in\mathcal{O}_{k+1}$,
the spectrum of $D\Phi_k(X,Y)$ cannot contain $0$ since 

$$\lvert D\Phi_k(X,Y)\rvert \geq 1-2\zeta_{k+1}.$$

\noindent Therefore $\Phi_k$ is injective on $\mathcal{O}_{k+1}$ and $\Phi_k^{-1}$ is defined on $\Phi_k(\mathcal{O}_{k+1})\subset \mathcal{O}_k$. Property \ref{DPhik} also comes from the estimate \eqref{DGk}.

\bigskip
\noindent
Moreover, $\Phi_k$ satisfies

$$\Phi_k^*(N_{k+1}+R_{k+1})=N_k+R_k$$

\noindent 
where 

$$N_{k+1}=N_k+T^{m_{k+1}-1} R_k-\bar{R}_{k}$$

\noindent and 

\begin{equation}\begin{split}
(I+DG_k)R_{k+1}&=-DG_k (T^{m_{k+1}-1} R_k-\bar{R}_{k}) + (R_k-T^{m_{k+1}-1} R_k) +DR_kG_k\\
& + \sum_{\lvert I\rvert \geq 2} \frac{1}{I_1!\dots I_{d+n}!} \partial_I (R_k+N_k) \cdot G_k^I
\end{split}\end{equation}

\noindent Thus, under the property \ref{mkcroiss} of Assumption \ref{gamma+omega}, $R_{k+1}$ has order $m_{k+1}$. Moreover, for any vector field $F$ of order $m_k$ and any index $I\in \mathbb{N}^{d+n}$,

\begin{equation}\begin{split}
&\lvert \lvert \partial_I F\rvert \rvert _{r_{k+1},\delta_{k+1}}\\
&=\sum_{P,\lvert Q\rvert \geq m_k}\left \lvert (iP_1)^{I_1}\dots (iP_d)^{I_d} \frac{Q_1!}{(Q_1-I_{d+1})!}\dots \frac{Q_n!}{(Q_n-I_{d+n})!} F_{P,Q} \right\lvert 
e^{\lvert P\rvert r_{k+1} } \delta_{k+1}^{\lvert Q\rvert -\lvert I \rvert}\\
&\leq (\delta_{k+1}/\delta_k)^{m_k}(r_k-r_{k+1})^{-i_1-\dots -i_d-d} (\ln \delta_k-\ln \delta_{k+1})^{-i_{d+1}-\dots -i_{d+n}-n} \lvert \lvert F\rvert \rvert _{r_k,\delta_k}\\
&\leq (\delta_{k+1}/\delta_k)^{m_k}(C(n,d)g(m_k))^{\lvert I\rvert +1+n}\lvert \lvert F\rvert \rvert _{r_k,\delta_k}\\
\end{split}
\end{equation}

\noindent which, applied to $N_k$ and $R_k$, implies

\begin{equation}\label{I+DGk}\begin{split}
&\lvert \lvert (I+DG_k)R_{k+1}\rvert \rvert_{r_{k+1},\delta_{k+1}} \\
&\leq \zeta_k\zeta_{k+1} +(\delta_{k+1}/\delta_k)^{m_{k+1}} \zeta_k + C(n,d) (\delta_{k+1}/\delta_k)^{m_k}g(m_k)^{n+2} \zeta_k \zeta_{k+1} \\
&+ (\delta_{k+1}/\delta_k)^{m_k}\sum_{\lvert I\rvert \geq 2} \frac{1}{I_1!\dots I_{d+n}!}(C(n,d)g(m_k))^{\lvert I\rvert +1+n}\zeta_{k+1}^{\lvert I\rvert } (\zeta_k+ \eta_k)\\
&\leq \zeta_k\zeta_{k+1} +\frac{1}{10} \zeta_{k+1} + \zeta_k \zeta_{k+1} \\
&+ (\delta_{k+1}/\delta_k)^{m_k}g(m_k)^{n+2}\sum_{\lvert I\rvert \geq 2} \frac{1}{I_1!\dots I_{d+n}!}(C(S,n,d)\zeta_{k})^{\lvert I\rvert-1 }\zeta_{k+1}\\
&\leq \frac{1}{2}\zeta_{k+1}
\end{split}\end{equation}

\bigskip
\noindent Now 

\begin{equation}\begin{split}
\lvert \lvert (I+DG_k)R_{k+1}\rvert \rvert_{r_{k+1},\delta_{k+1}} &\geq \lvert \lvert R_{k+1}\rvert \rvert_{r_{k+1},\delta_{k+1}} - 
\lvert \lvert DG_kR_{k+1}\rvert \rvert_{r_{k+1},\delta_{k+1}} \\
&\geq (1-\zeta_{k+1})\lvert \lvert R_{k+1}\rvert \rvert_{r_{k+1},\delta_{k+1}}
\end{split}\end{equation}

\noindent which, together with \eqref{I+DGk}, finally implies that $\lvert \lvert R_{k+1}\rvert \rvert_{r_{k+1},\delta_{k+1}} \leq \zeta_{k+1}$, whence property \ref{2}.

\bigskip
\noindent The vector field $N_{k+1}$ satisfies the following estimate:

\begin{equation}
\lvert \lvert N_{k+1}-S\rvert \rvert_{r_{k+1},\delta_{k+1}}\leq \lvert\lvert  N_{k}-S\rvert \rvert_{r_{k},\delta_{k}}+ \lvert\lvert  {R}_k\rvert  \rvert_{r_k,\delta_k} \leq \eta_k +\zeta_k \leq \eta_{k+1}
\end{equation}

\noindent whence property \ref{1}. This concludes the proof of Proposition \ref{iteration}. $\Box$

\nocite{*}
\bibliographystyle{cdraifplain}
\bibliography{xampl}

\end{document}